\documentclass[12pt]{elsart}
\usepackage{amsmath,amsfonts,amssymb,amsbsy}
\usepackage[pdftex]{graphicx}

\begin{document}

\begin{frontmatter}
\title{ 
Mathematical modeling of linear viscoelastic impact: 
Application to drop impact testing of articular cartilage}
\author{I.I.~Argatov}
\ead{iva1@aber.ac.uk}
\address{Institute of Mathematics and Physics, Aberystwyth University, Ceredigion SY23 3BZ, Wales, UK}

\begin{abstract}
In recent years, a number of experimental studies have been conducted to investigate the mechanical behavior and damage mechanisms of articular cartilage under impact loading. Some experimentally observed results have been explained using a non-linear viscoelastic impact model. At the same time, there is the need of simple mathematical models, which allow comparing experimental results obtained in drop impact testing with impact loads of different weights and incident velocities. The objective of this study was to investigate theoretically whether the main features of articular impact could be qualitatively predicted using a linear viscoelastic theory or the linear biphasic theory. In the present paper, exact analytical solutions are obtained for the main parameters of the Kelvin--Voigt and Maxwell impact models. Perturbation analysis of the impact process according to the standard viscoelastic solid model is performed. Asymptotic solutions are obtained for the drop weight impact test. The dependence of the coefficient of restitution on the impactor parameters has been studied in detail. 

\end{abstract}

\begin{keyword}
Impact contact problem \sep blunt indenter \sep asymptotic model \sep coefficient of restitution
\end{keyword}
\end{frontmatter}

\newpage
\section*{Nomenclature}

\begin{tabular}{rl}
$b$ & damper coefficient \\
$D$ & discriminant of the characteristic equation  \\
$E_{\rm dyn}$ & incremental dynamic modulus  \\
$E_{\rm max}$ & maximum incremental dynamic modulus  \\
$E_{10}$ & modulus at stresses of 10~MPa  \\
$e_*$ & coefficient of restitution \\
$F$ &  contact force  \\
$F_M$ & maximum contact force \\
$g$ &  gravitational acceleration  \\
$h$ &  cartilage layer thickness  \\
$h_0$ &  drop height of the impactor  \\
$H_A$ & aggregate modulus \\
$k$ & stiffness coefficient \\
$k_1$, $k_2$ & spring stiffnesses in the standard solid model \\
$k_0$ & instantaneous stiffness  \\
$k_\infty$ & long-term stiffness  \\
$m$ & impactor mass \\
$t$ & time variable \\
$t_c$ & impact duration \\
$t_m$ & time to maximum displacement \\
$t_M$ & time to maximum contact force \\
$v_0$ &   initial impact velocity  \\
$x$ &   displacement  \\
$\dot{x}$ &   velocity  \\
$\ddot{x}$ &  acceleration  \\
$x_m$ & maximum displacement
\end{tabular}
\vfill
$${}$$

\begin{tabular}{rl}
$\beta$ &  damping coefficient in the Kelvin--Voigt model  \\
$\beta_1$ & real part of complex roots of the characteristic equation  \\
$\Delta m$ & percentage increase in mass of cartilage sample  \\
$\epsilon$ & strain  \\
$\varepsilon_0$ &  non-dimensional parameter accounting for the gravitational effect \\
$\zeta$ &  loss factor in the Maxwell model  \\
$\zeta_1$ & imaginary part of complex roots of the characteristic equation  \\
$\eta$ &  loss factor in the Kelvin--Voigt model  \\
$\kappa$ & cartilage permeability  \\
$\varkappa_1$, $\varkappa_2$ & spring stiffnesses in the standard solid model \\
$\lambda$ &  Lam\'e coefficient  \\
$\lambda_1$ &  root of the characteristic equation  \\
$\Lambda$ &  non-dimensional parameter in the standard solid model  \\
$\mu$ &  Lam\'e coefficient  \\
$\xi$ &  non-dimensional displacement  \\
$\rho$ &  ratio of the long-term and instantaneous stiffnesses   \\
$\sigma$ &  stress   \\
$\tau$ &  non-dimensional time  \\
$\tau_D$ &  typical diffusion time  \\
$\tau_R$ &  relaxation time  \\
$\Psi(\tau)$ & dimensionless relaxation function  \\
$\omega$ &  angular frequency of damped oscillations  \\
$\omega_0$ &  angular frequency of undamped oscillations  \\
\end{tabular}

\setcounter{equation}{0}
\section{Introduction} 
\label{1dsSectionI}

Articular cartilage is a soft hydrated tissue covering the end of each bone at the joints. Cartilage has no known function other than maintaining mechanical competence of joints, allowing bones to move against one another without friction. But there is no need to underline its significance to health of a human body, since almost all the load transmitted by a human joint goes through the articular cartilage, and it prevents biomechanical damage caused by severe loading including impact loading. It is believed that severe articular impact can initiate post-traumatic arthritis \cite{JeffreyGregoryAspden1995,QuinnAllen2001}. An impact loading of the joint constitutes the action of extremely high non-physiological loads applied very rapidly (for instance, due to a car accident, sports injury, or a fall from a height). 

\begin{figure}[h!]
    \centering
    \includegraphics [scale=0.35]{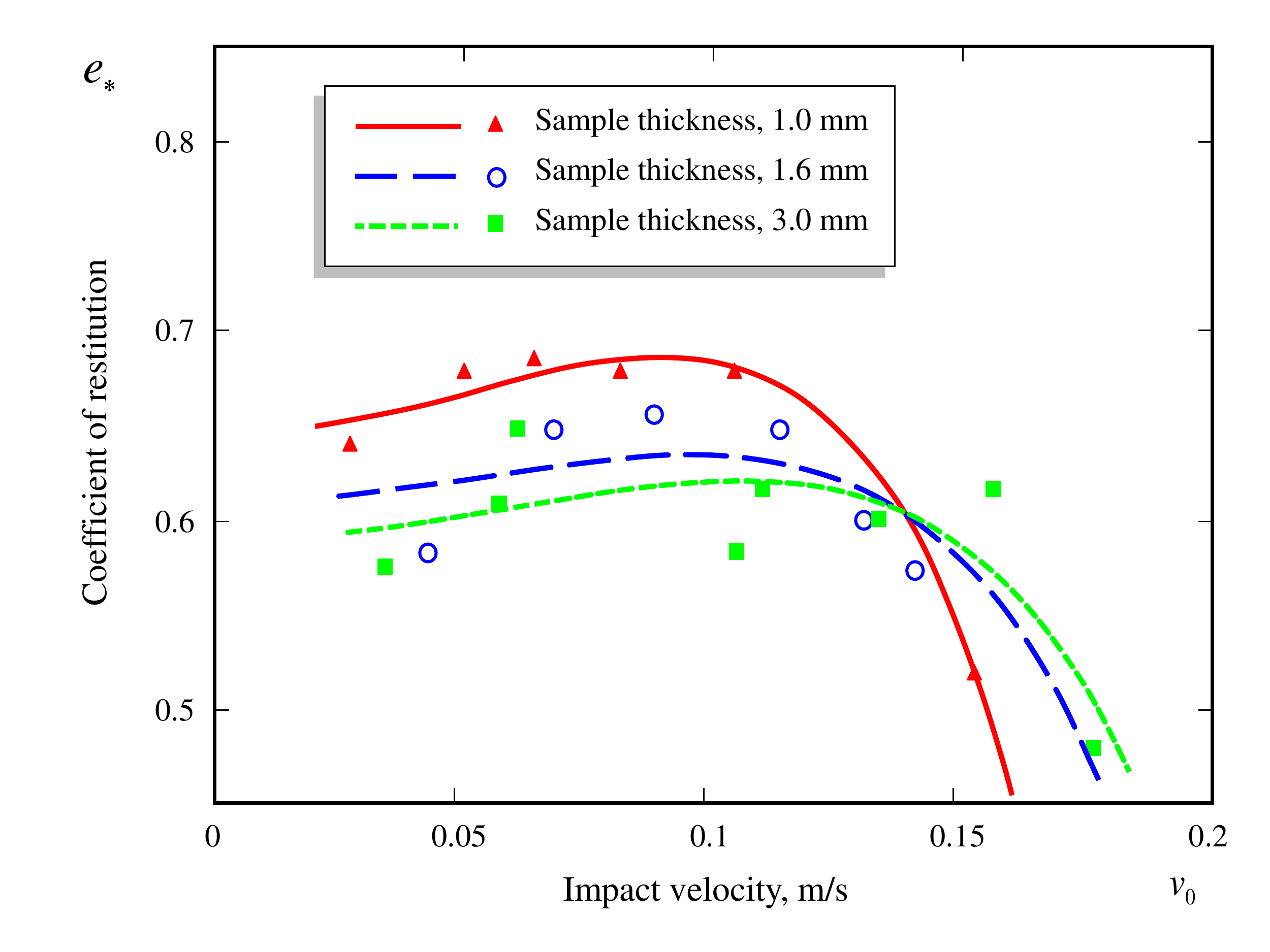}
    \caption{Coefficient of restitution $e_*$ versus the impact velocity $v_0$ for articular cartilage samples of different thicknesses. Based on the experimental data and fitting curves obtained in \cite{Varga2007}.
    }
    \label{Varga2007.pdf}
\end{figure}

In recent years, a number of experimental studies have been conducted to investigate the mechanical behavior and damage mechanisms of articular cartilage under impact loading 
\cite{AtkinsonHautAltiero1998,VerteramoSeedhom2007,BurginAspden2008}. In particular, the experimental data on relative dissipation of the impact energy $\Delta E/E_0$ versus overall impactor energy $E_0$ obtained in \cite{Varga2007} were fitted with quadratic curves. Here, $E_0=mv_0^2/2$, 
$\Delta E=m(v_1^2-v_0^2)/2$, $v_0$ and $v_1$ are the initial impact and rebound velocities, respectively, $m$ is the impactor mass. Since, $v_1=-e_* v_0$, where $e_*$ is the coefficient of restitution, we easily get $\Delta E/E_0=1-e_*^2$. Thus, the experimental data and fitting curves for dissipation of the impact energy \cite{Varga2007} can be recalculated in terms of the coefficient of restitution as presented in 
Fig.~\ref{Varga2007.pdf}, which shows a non-monotonic dependence of $e_*$ on $v_0$. 
Some experimentally observed results have been explained using a non-linear viscoelastic impact model \cite{Edelsten2010}. At the same time, there is the need of a simple mathematical model, which allows comparing experimental results obtained in drop impact testing with impact loads of different weights and incident velocities. 

A variety of mathematical models were suggested to describe the stress-strain response of articular cartilage that represents a multiphasic, structurally complex material possessing viscoelastic properties.
It is long known that articular cartilage possesses viscoelastic properties \cite{HayesMockros1971,Lau_et_al_2008}, though there is no direct correspondence between viscoelastic parameters and parameters of the biphasic/poroelastic models of cartilage.
The biphasic theory \cite{MowKueiLaiArmstrong1980}, which models the tissue as a mixture of a solid phase and a fluid phase, has demonstrated very good agreement with experimental results in the creep and stress relaxation tests \cite{SoltzAteshian2000}. The objective of this study was to investigate theoretically whether the main features of articular impact observed in \cite{Varga2007,Edelsten2010} could be qualitatively predicted using a linear viscoelastic theory or the linear biphasic theory.

The rest of the paper is organized as follows. In Sections~\ref{1dsSection1} and \ref{1dsSection2}, we consider in detail the viscoelastic Kelvin--Voigt and Maxwell impact models, respectively. Since some elements of the presented solutions are known in the literature, we pay a particular attention to the evaluation of the contact force, $F(t)$, and impactor displacement, $x(t)$, at the time moments $t_M$ and $t_m$, when the force and displacement reach their maxima, $F_M$ and $x_m$, respectively. 
In Section~\ref{1dsSection3}, we outline a closed form solution of the impact equation in the case of standard solid model. In order to get analytical approximations, we consider the standard solid model as a perturbation of the Kelvin--Voigt (Section~\ref{1dsSection4}) or the Maxwell model (Section~\ref{1dsSection5}). In particular, simple analytical approximations are derived for the impact duration, $t_c$, and for the coefficient of restitution, $e_*$. 
In Sections~\ref{1dsSection05} and \ref{1dsSection06}, we consider the influence of the gravity effect on these parameters in the framework of the Kelvin--Voigt and Maxwell models for drop weight impact. 
In Section~\ref{1dsSection07}, we develop an asymptotic model for the force-displacement relationship in the indentation problem for a thin biphasic layer corresponding to the conditions of the so-called blunt impact, when the specimen thickness is much smaller than the radius of a flat-ended cylindrical impactor.
An example of application of the developed linear theory of viscoelastic impact for analyzing experimental data is given in Section~\ref{1dsSection10}.
Finally, in Sections~\ref{1dsSectionD} and \ref{1dsSectionC}, we outline a discussion of the results obtained and formulate our conclusions.

\section{Viscoelastic Kelvin--Voigt impact model}
\label{1dsSection1}

\begin{figure}[h!]
    \centering
    \includegraphics [scale=0.35]{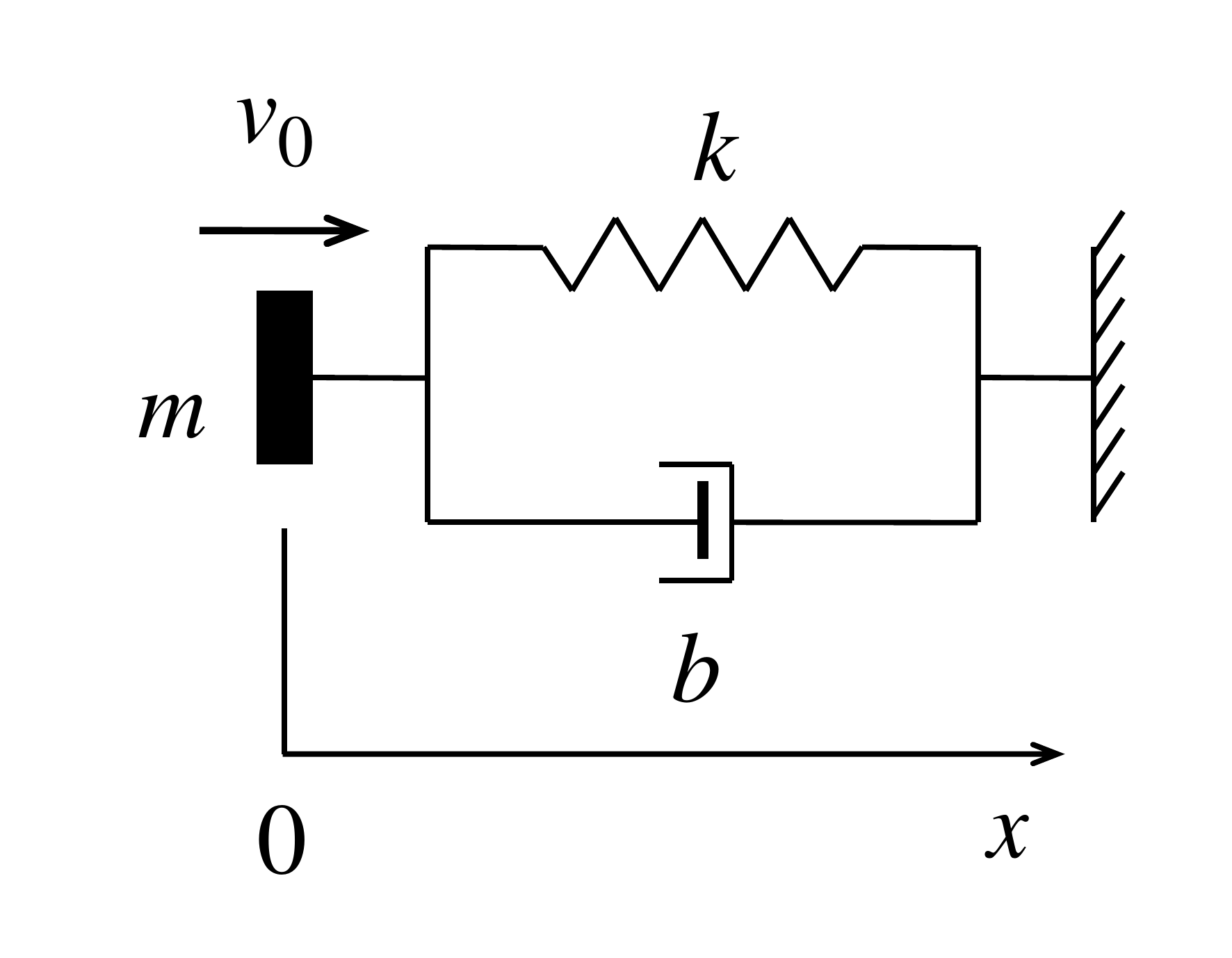}
    \caption{Impact viscoelastic Kelvin--Voigt model.}
    \label{Kelvin-Voigt_Impact.pdf}
\end{figure}

In this section, the deformation of articular cartilage layer is modeled schematically as a parallel combination of linear spring $k$ and dashpot $b$ (Fig.~\ref{Kelvin-Voigt_Impact.pdf}). Dynamic balance between the force of cartilage reaction
\begin{equation}
F=kx+b\dot{x}
\label{1vI(1.0)}
\end{equation}
and the force of body inertia $m\ddot{x}$ governs the development of collision. According to Newton's second law, the differential equation of the impact has the form
\begin{equation}
m\ddot{x}+b\dot{x}+kx=0, \quad t\in[0,t_c],
\label{1vI(1.1)}
\end{equation}
where $t_c$ is the contact duration, that is $t_c$ denotes the instant, when the cartilage reaction force changes its sign, or the indenter acceleration vanishes. 

The initial conditions for Eq.~(\ref{1vI(1.1)}) are as follows:
\begin{equation}
x(0)=0, \quad \dot{x}(0)=v_0.
\label{1vI(1.2)}
\end{equation}

The impact duration is determined by the condition 
\begin{equation}
kx+b\dot{x}\bigr\vert_{t=t_c}=0,
\label{1vI(1.3)}
\end{equation}
or, in view of Eq.~(\ref{1vI(1.1)}), by the condition
\begin{equation}
\ddot{x}\bigr\vert_{t=t_c}=0.
\label{1vI(1.4)}
\end{equation}

The impact problem (\ref{1vI(1.1)}), (\ref{1vI(1.2)}) has the following well-known solution \cite{WinemanRajagopal2000}:
\begin{equation}
x(t)=\frac{v_0}{\omega}e^{-\beta t}\sin\omega t, \quad t\in[0,t_c].
\label{1vI(1.5)}
\end{equation}
Here we used the notation
\begin{equation}
\omega_0^2=\frac{k}{m},\quad \omega^2=\omega_0^2-\beta^2, \quad \beta=\frac{b}{2m}.
\label{1vI(1.6)}
\end{equation}
We assume that $\omega_0>\beta$.

\begin{figure}[h!]
    \centering
    \includegraphics [scale=0.35]{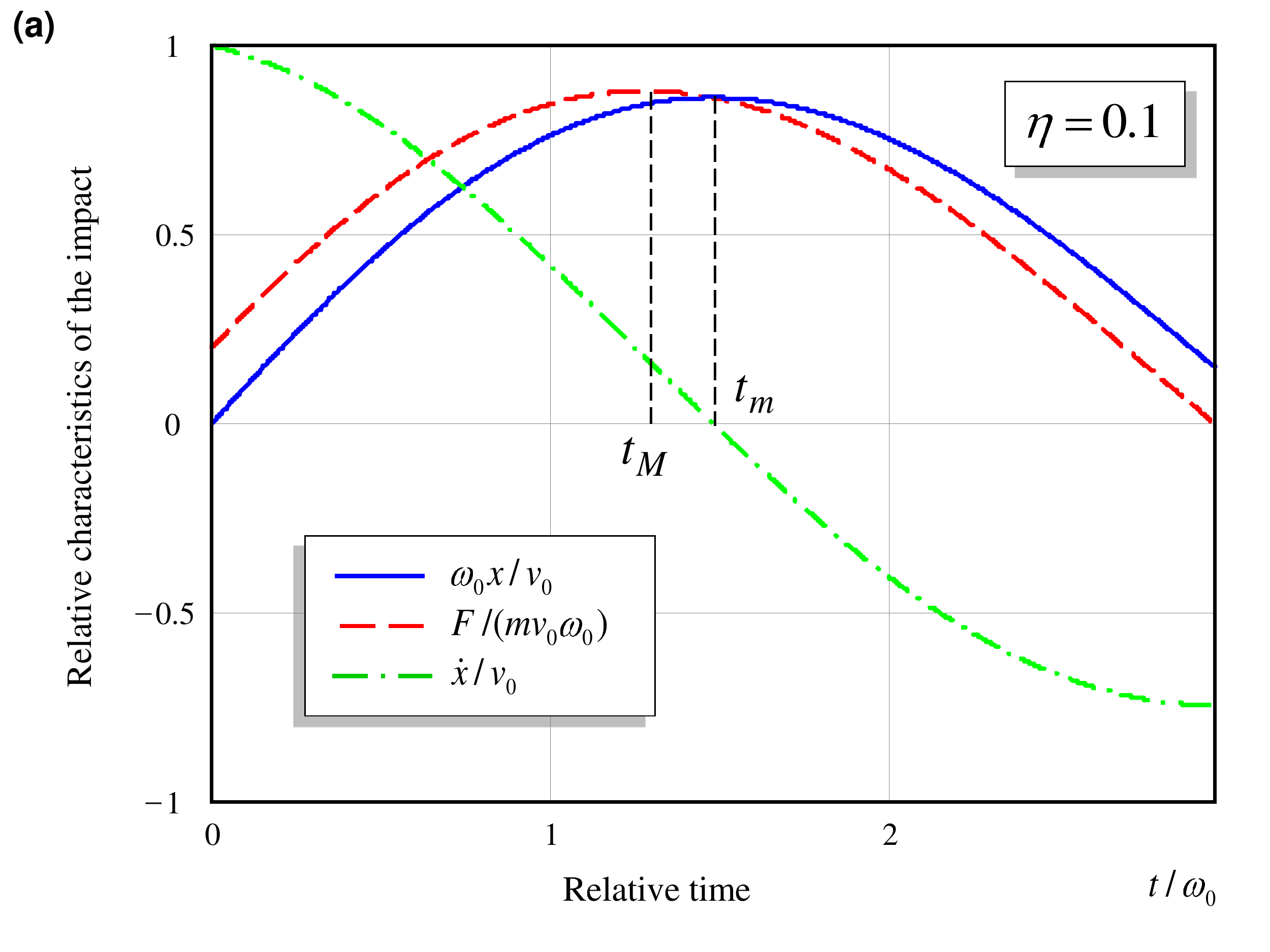}
    \includegraphics [scale=0.35]{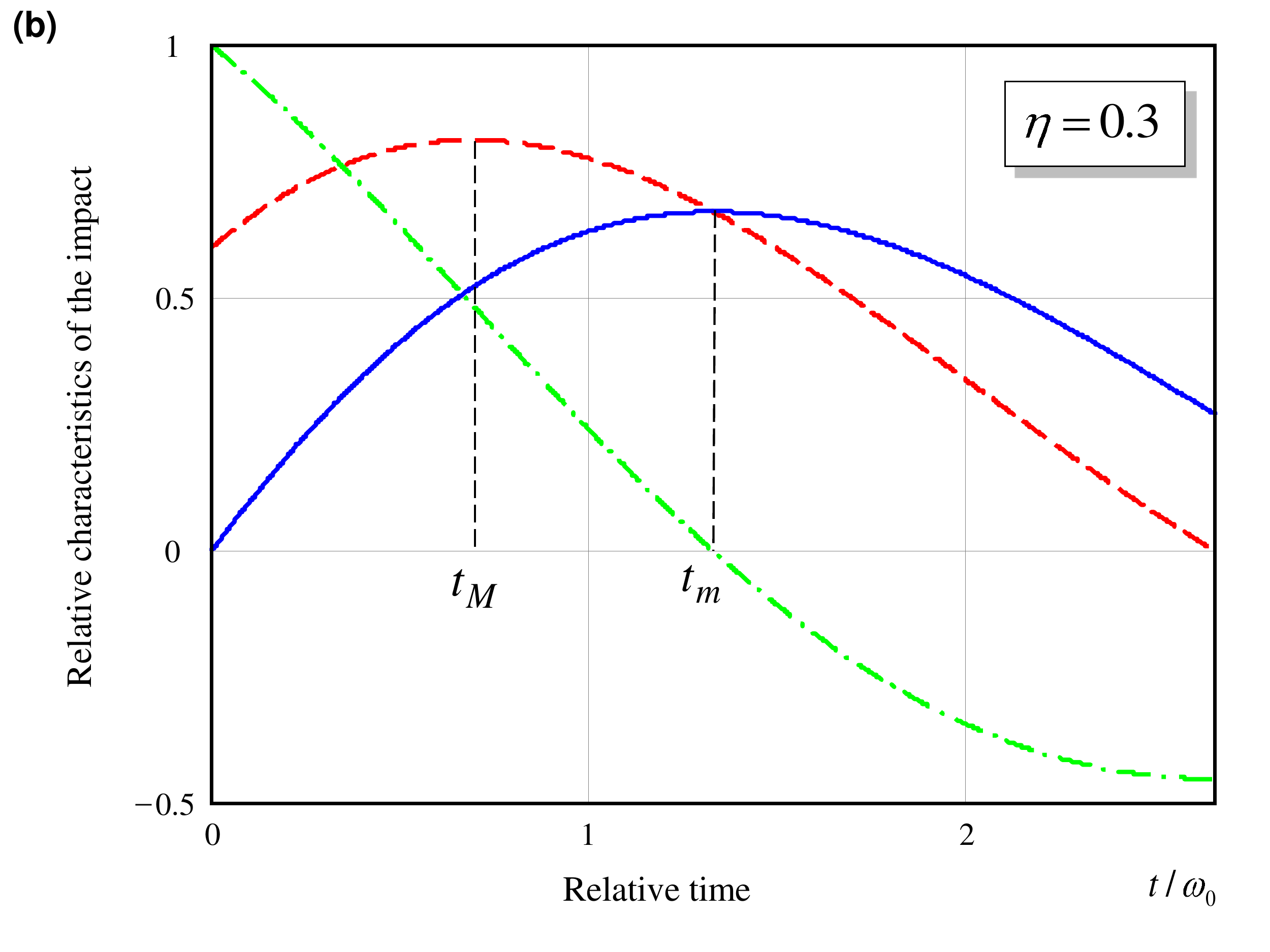}
    \includegraphics [scale=0.35]{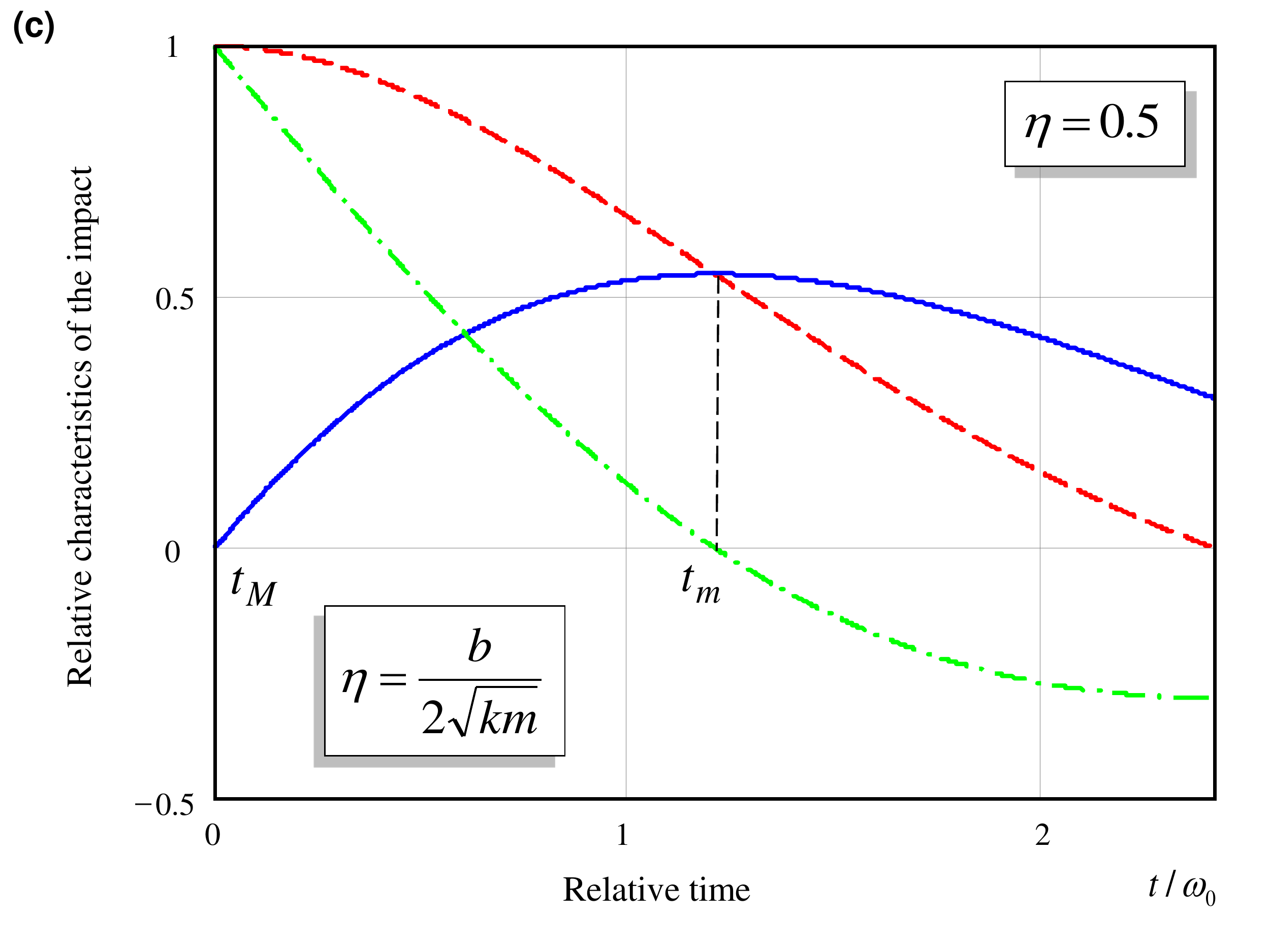}
    \caption{Viscoelastic Kelvin--Voigt impact model. Behavior of the main impact variables with time for the following values of the damping ratio: $\eta=0{.}1$ (a), $\eta=0{.}3$ (b), $\eta=0{.}5$ (c).}
    \label{Kelvin-Voigt_exact.pdf}
\end{figure}

Fig.~\ref{Kelvin-Voigt_exact.pdf} shows the behavior of the dimensionless quantities 
$\omega_0 x/v_0$, $F/(mv_0\omega_0)$, and $\dot{x}/v_0$ with respect to time. Observe that the time moment $t_M$, when the contact force reaches its maximum, approaches the initial moment of impact as the damping ratio $\eta$ increases. 

By differentiating the both sides of Eq.~(\ref{1vI(1.5)}) with respect to $t$, we get
\begin{equation}
\dot{x}(t)=\frac{v_0}{\omega}e^{-\beta t}\bigl(
\omega\cos\omega t-\beta\sin\omega t\bigr),
\label{1vI(1.7a)}
\end{equation}
\begin{equation}
\ddot{x}(t)=-\frac{v_0}{\omega}e^{-\beta t}\bigl[
(\omega^2-\beta^2)\sin\omega t+2\beta\omega\cos\omega t\bigr].
\label{1vI(1.7b)}
\end{equation}

After solving Eq.~(\ref{1vI(1.4)}) for $t_c$ in view of (\ref{1vI(1.7b)}), the following expression for the impact duration can be obtained \cite{ButcherSegalman2000}:
\begin{equation}
t_c=\frac{1}{\omega}{\rm Atan\,}\frac{-2\beta\omega}{\omega^2-\beta^2},
\label{1vI(1.8)}
\end{equation}
where the first positive value of the many-valued ${\rm Atan}$ function should be taken. 

Using properties of the ${\rm Atan}$ function, we rewrite Eq.~(\ref{1vI(1.8)}) as follows \cite{Popov2010}:
\begin{equation}
t_c=\frac{1}{\omega}
\left\{
\begin{array}{l}
\displaystyle
\pi-{\,\rm atan\,}\frac{2\beta\omega}{\omega^2-\beta^2},\quad \beta<\omega, \\
\displaystyle
{\,\rm atan\,}\frac{2\beta\omega}{\omega^2-\beta^2}, \quad \omega>\beta.
\end{array}
\right.
\label{1vI(1.9)}
\end{equation}
Here, ${\,\rm atan\,}(z)$ is the principal branch of the arctangent function ${\,\rm Atan\,}(z)$.

Finally, using properties of the ${\,\rm atan\,}$ function, we can rewrite formula (\ref{1vI(1.9)}) in a more simple form as 
\begin{equation}
t_c=\frac{2}{\omega}{\,\rm atan\,}\frac{\omega}{\beta}.
\label{1vI(1.9a)}
\end{equation}

Let $\eta$ denote the loss factor, i.\,e.,
\begin{equation}
\eta=\frac{\beta}{\omega_0}.
\label{1vI(1.10a)}
\end{equation}
Then, Eq.~(\ref{1vI(1.9a)}) can be rewritten as
\begin{equation}
t_c=\frac{2}{\omega_0\sqrt{1-\eta^2}}{\,\rm atan\,}\frac{\sqrt{1-\eta^2}}{\eta}.
\label{1vI(1.10)}
\end{equation}
Recall that we assume that $\eta\in[0,1]$. Also, note that in view of the notation (\ref{1vI(1.6)}), we have
\begin{equation}
\eta=\frac{b}{2\sqrt{km}}.
\label{1vI(1.10b)}
\end{equation}

The velocity of the indenter at separation can be obtained by the substitution of (\ref{1vI(1.10)}) into (\ref{1vI(1.7a)}) in the following form:
\begin{eqnarray}
\dot{x}(t_c) & = & -v_0\exp(-\beta t_c) 
\label{1vI(1.11)} \\
{} & = & -v_0\exp\biggl\{- \frac{2\eta}{\sqrt{1-\eta^2}}{\,\rm atan\,}\frac{\sqrt{1-\eta^2}}{\eta}\biggr\}.
\label{1vI(1.12)}
\end{eqnarray}

From Eq.~(\ref{1vI(1.11)}), it follows that the coefficient of restitution, $e_*$, which is defined as the ratio of the velocity at separation $\vert \dot{x}(t_c)\vert$ to the velocity of the indenter at incidence 
$\vert \dot{x}(0)\vert=v_0$, is given by
\begin{eqnarray}
e_* & = & \exp\Bigl\{- \frac{2\beta}{\omega}{\,\rm atan\,}\frac{\omega}{\beta}\Bigr\}
\label{1vI(1.13)} \\
{} & = & \exp\biggl\{- \frac{2\eta}{\sqrt{1-\eta^2}}{\,\rm atan\,}\frac{\sqrt{1-\eta^2}}{\eta}\biggr\}.
\label{1vI(1.14)}
\end{eqnarray}

The peak value of the indenter penetration occurs at the instant $t=t_m$, when $\dot{x}(t_m)=0$. In view of 
(\ref{1vI(1.7a)}), we have
\begin{eqnarray}
t_m & = & \frac{1}{\omega}{\,\rm atan\,}\frac{\omega}{\beta}
\label{1vI(1.15)} \\
{} & = & \frac{1}{\omega_0\sqrt{1-\eta^2}}{\,\rm arcsin\,}\sqrt{1-\eta^2}.
\label{1vI(1.16)}
\end{eqnarray}
Substituting the expression (\ref{1vI(1.15)}) into Eq.~(\ref{1vI(1.5)}), we obtain the maximum penetration 
$x_m=x(t_m)$ in the form
\begin{eqnarray}
x_m & = & \frac{v_0}{\omega_0}\exp\Bigl(-\frac{\beta}{\omega}{\,\rm atan\,}\frac{\omega}{\beta}\Bigr)
\label{1vI(1.17)} \\
{} & = & \frac{v_0}{\omega_0}\exp\Bigl(-\frac{\eta}{\sqrt{1-\eta^2}}{\,\rm arcsin\,}\sqrt{1-\eta^2}
\Bigr).
\label{1vI(1.18)}
\end{eqnarray}

Note that from (\ref{1vI(1.10)}) and (\ref{1vI(1.16)}), it is readily seen that $t_m=t_c/2$.

The peak value of the contact force $F$ occurs at the instant $t=t_M$, when $\dot{F}(t_M)=0$. According to Eqs.~(\ref{1vI(1.5)}), (\ref{1vI(1.7a)}), we obtain
\begin{eqnarray}
F(t) & = & \frac{mv_0}{\omega}\exp(-\beta t)\bigl[(\omega^2-\beta^2)\sin\omega t+2\beta\omega\cos\omega t\bigr]
\nonumber \\
{} & = & m v_0\omega_0 \exp(-\eta\omega_0 t)\biggl(\frac{(1-2\eta^2)}{\sqrt{1-\eta^2}}
\sin\omega_0\sqrt{1-\eta^2}t
+2\eta\cos\omega_0\sqrt{1-\eta^2}t
\biggr).
\label{1vI(1.19)}
\end{eqnarray}

Differentiating the previous expression, we can reduce the equation $\dot{F}(t_M)=0$ to the following one:
$$\sqrt{1-\eta^2}(1-4\eta^2)\cos(t_M\omega_0\sqrt{1-\eta^2})+
\eta(4\eta^2-3)\sin(t_M\omega_0\sqrt{1-\eta^2})=0.$$

Thus, for $\eta\in(0,0{.}5)$, we obtain 
\begin{eqnarray}
t_M & = & \frac{1}{\omega}{\,\rm atan\,}\frac{\omega(\omega^2-3\beta^2)}{\beta
(3\omega^2-\beta^2)}
\nonumber \\
{} & = & \frac{1}{\omega_0\sqrt{1-\eta^2}}{\,\rm atan\,}\frac{\sqrt{1-\eta^2}(1-4\eta^2)}{\eta(3-4\eta^2)}.
\label{1vI(1.20)}
\end{eqnarray}
For $\eta\in(0{.}5,1)$, the maximum value of the contact force $F_M=F(t_M)$ takes place at the initial instant $t=0$.

Substituting (\ref{1vI(1.20)}) into Eq.~(\ref{1vI(1.19)}), we get
\begin{eqnarray}
F_M & = & m v_0\omega_0 \exp\biggl(-\frac{\eta}{\sqrt{1-\eta^2}}
{\,\rm atan\,}\frac{\sqrt{1-\eta^2}(1-4\eta^2)}{\eta(3-4\eta^2)}\biggr),\quad \eta\in(0,0{.}5), \\
\label{1vI(1.21)}
F_M & = & m v_0\omega_0 2\eta,\quad \eta\in(0{.}5,1).
\label{1vI(1.22)}
\end{eqnarray}

Note that the function $F_M(\eta)$ defined by Eqs.~(\ref{1vI(1.21)}) and (\ref{1vI(1.22)}) is continuously differentiable. 

\begin{figure}[h!]
    \centering
    \includegraphics [scale=0.35]{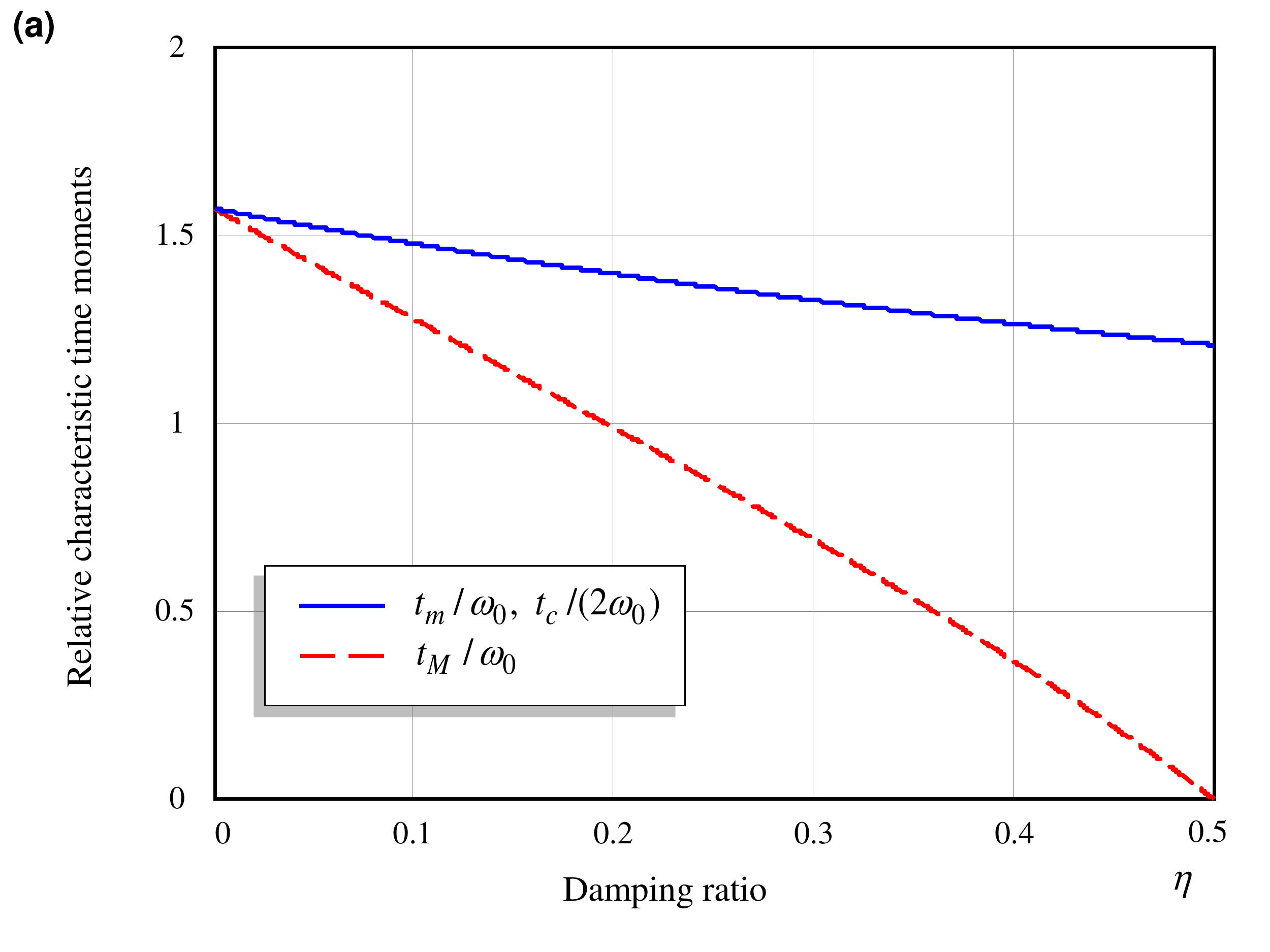}
    \includegraphics [scale=0.35]{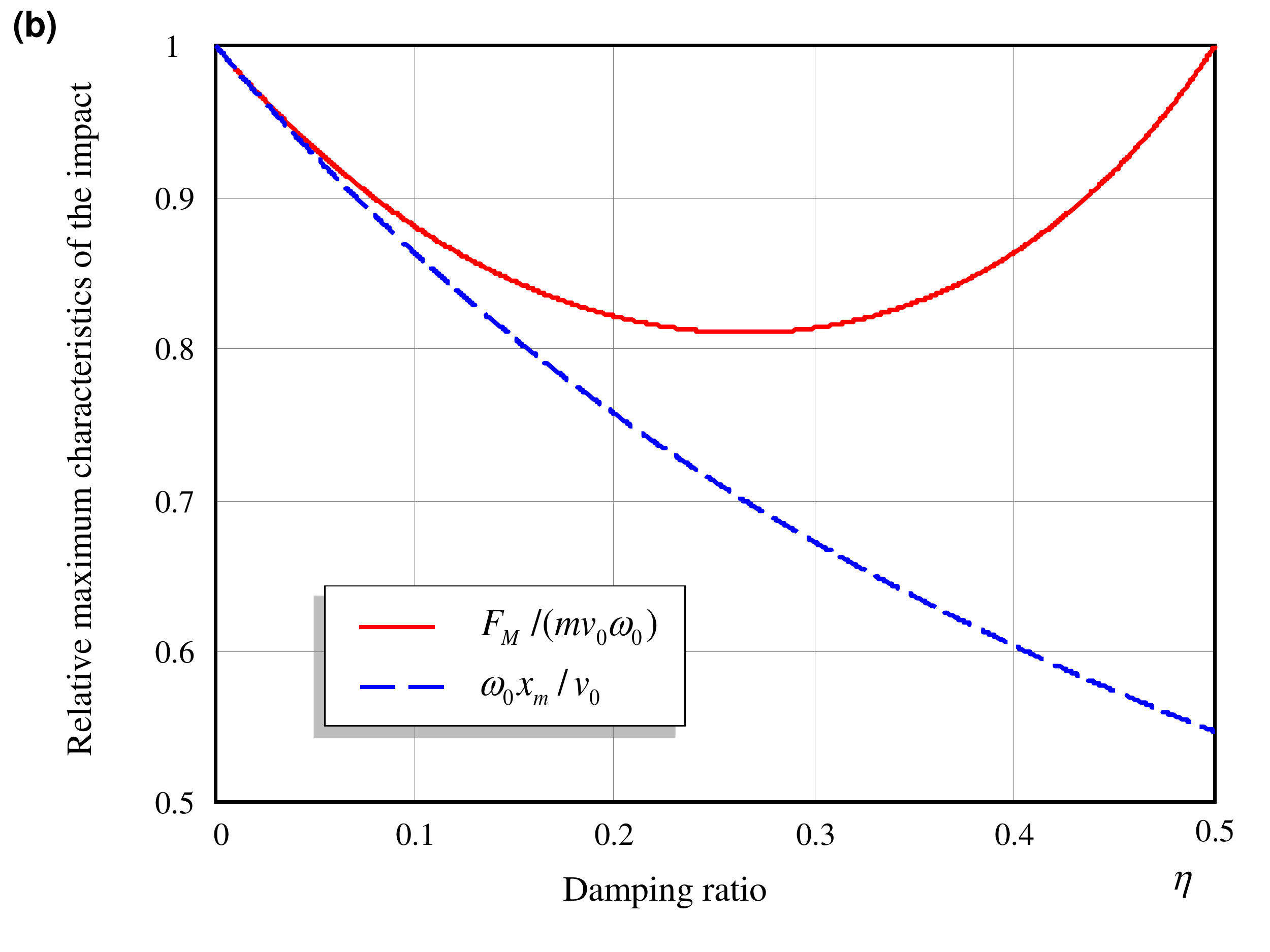}
    \caption{Viscoelastic Kelvin--Voigt impact model. Behavior of the main impact parameters 
    $t_m$, $t_c$, $t_M$ (a) and $x_m$, $F_M$ (b) with the damping ratio.}
    \label{Kelvin-Voigt_tm+tM.pdf}
\end{figure}

Fig.~\ref{Kelvin-Voigt_tm+tM.pdf}a shows the monotonic behavior of the dimensionless characteristic time moments $t_m/\omega_0$, $t_c/(2\omega_0)$, $t_M/\omega_0$ with the damping ratio $\eta$. Recall that $t_m=t_c/2$. The variations of the relative maximum contact force $F_M/(mv_0\omega_0)$ and displacement $\omega_0 x_m/v_0$ are presented in
Fig.~\ref{Kelvin-Voigt_tm+tM.pdf}b. It is interesting to observe the non monotonic behavior of $F_M$ with the minimum at $\eta\approx 0{.}26493$.

\section{Viscoelastic Maxwell impact model}
\label{1dsSection2}

\begin{figure}[h!]
    \centering
    \includegraphics [scale=0.35]{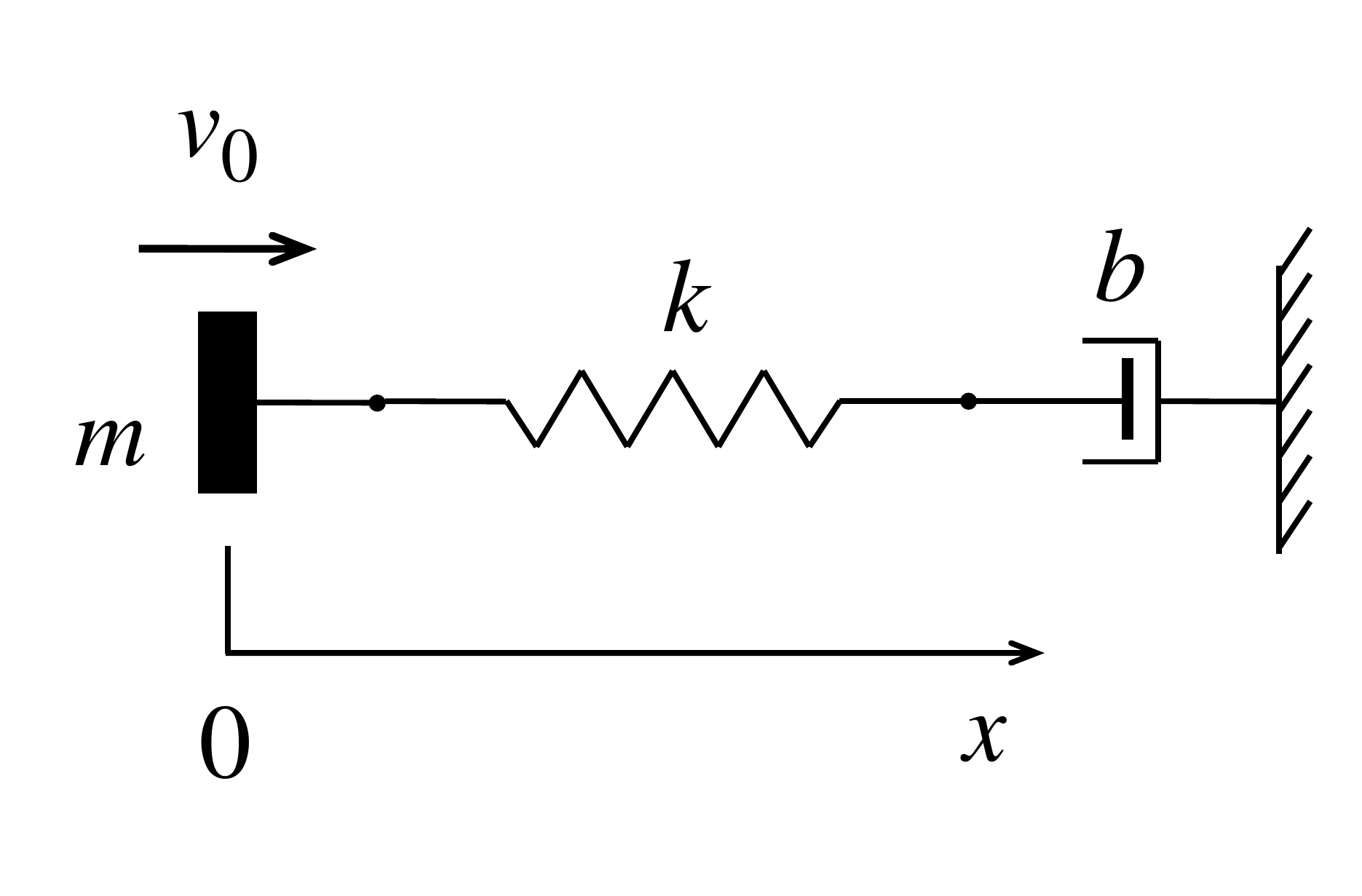}
    \caption{Impact viscoelastic Maxwell model.}
    \label{Maxwell_Impact.pdf}
\end{figure}

Assuming that the cartilage layer's response to impact loading is modeled schematically as a serial combination of linear spring $k$ and dashpot $b$ (Fig.~\ref{Maxwell_Impact.pdf}). The force-displacement relation is given by the following differential equation \cite{WinemanRajagopal2000}:
\begin{equation}
\frac{\dot{F}}{k}+\frac{F}{b}=\dot{x}.
\label{1vI(2.1)}
\end{equation}

From (\ref{1vI(2.1)}), it follows that 
\begin{equation}
F=k\int\limits_0^t\exp\Bigl\{-\frac{k}{b}(t-\tau)\Bigr\}\frac{dx}{d\tau}(\tau)\,d\tau.
\label{1vI(2.2)}
\end{equation}

The differential equation of the impact $m\ddot{x}+F=0$ in view of (\ref{1vI(2.1)}) results in the third-order equation
\begin{equation}
\dddot{x}+\frac{k}{b}\ddot{x}+\frac{k}{m}\dot{x}=0
\label{1vI(2.3)}
\end{equation}
with the initial conditions
\begin{equation}
x(0)=0, \quad \dot{x}(0)=v_0,\quad \ddot{x}(0)=0.
\label{1vI(2.4)}
\end{equation}

The impact problem (\ref{1vI(2.3)}), (\ref{1vI(2.4)}) has the following solution \cite{ButcherSegalman2000,Stronge2000}:
\begin{equation}
x(t)=\frac{v_0}{\omega_0}\exp(-\zeta\omega_0 t)\Bigl\{\frac{\omega_0(1-2\zeta^2)}{\omega}\sin\omega t
-2\zeta\cos\omega t\Bigr\}+\frac{2\zeta v_0}{\omega_0},
\label{1vI(2.5)}
\end{equation}
\begin{equation}
\dot{x}(t)=v_0\exp(-\zeta\omega_0 t)\Bigl\{
\cos\omega t+\frac{\zeta\omega_0}{\omega}\sin\omega t\Bigr\}.
\label{1vI(2.5a)}
\end{equation}
Here we used the notation
\begin{equation}
\omega_0^2=\frac{k}{m},\quad \omega=\omega_0\sqrt{1-\zeta^2}, \quad \zeta=\frac{k}{2\omega_0 b}.
\label{1vI(2.6)}
\end{equation}

The variation of the contact force during the impact interaction is 
\begin{equation}
F=\frac{kv_0}{\omega}\exp(-\zeta\omega_0 t)\sin\omega t.
\label{1vI(2.7)}
\end{equation}

The impact duration $t_c$ is determined by the condition 
$F\bigr\vert_{t=t_c}=0$. Thus, according to (\ref{1vI(2.7)}), the following relation takes place \cite{ButcherSegalman2000,Stronge2000}:
\begin{equation}
t_c=\frac{\pi}{\omega}.
\label{1vI(2.8)}
\end{equation}

Substituting the value (\ref{1vI(2.8)}) into Eq.~(\ref{1vI(2.5a)}), one gets the coefficient of restitution in the form
\begin{equation}
e_*=\exp\Bigl(-\frac{\pi\zeta}{\sqrt{1-\zeta^2}}\Bigr).
\label{1vI(2.9)}
\end{equation}

\begin{figure}[h!]
    \centering
    \includegraphics [scale=0.35]{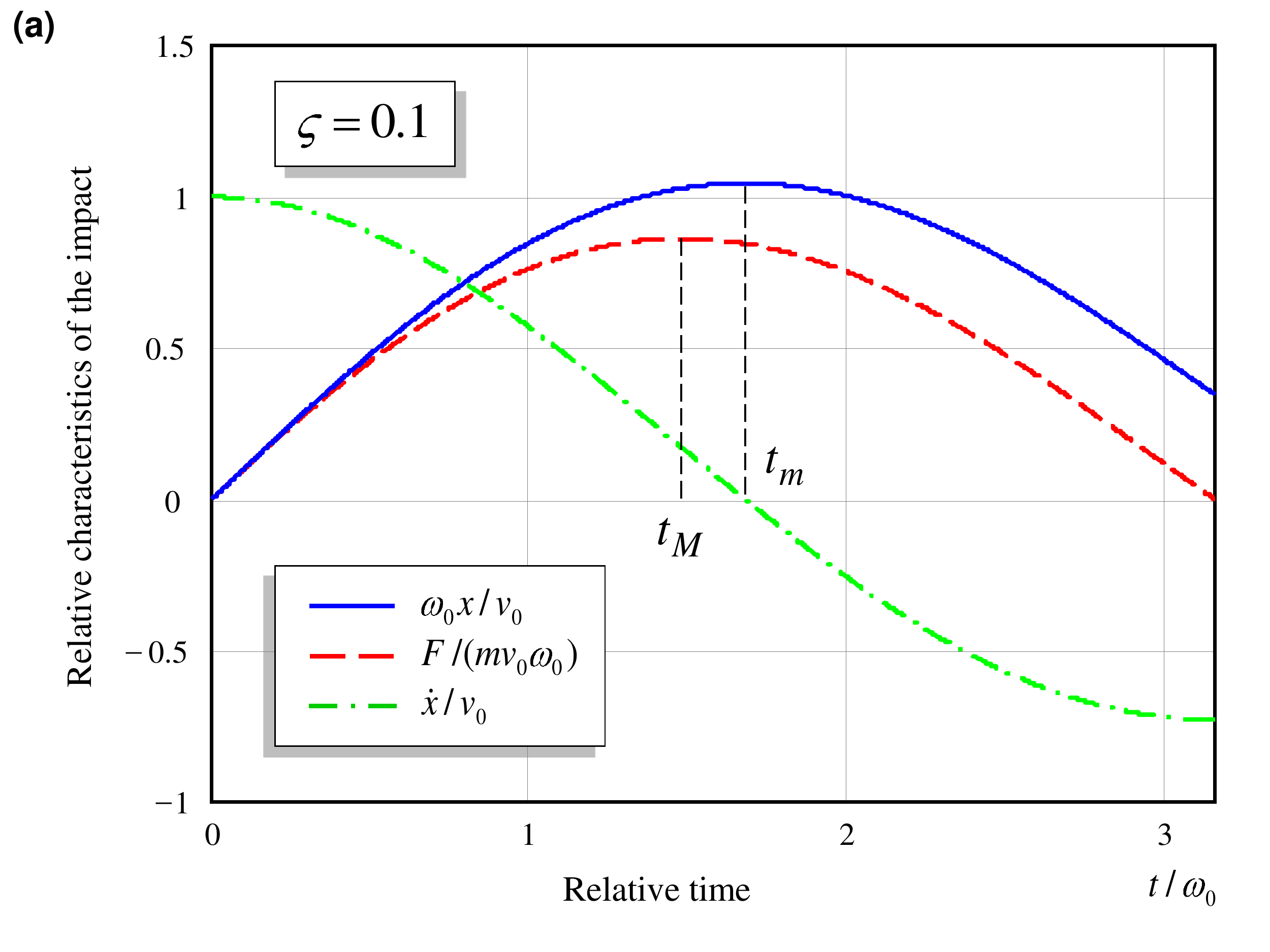}
    \includegraphics [scale=0.35]{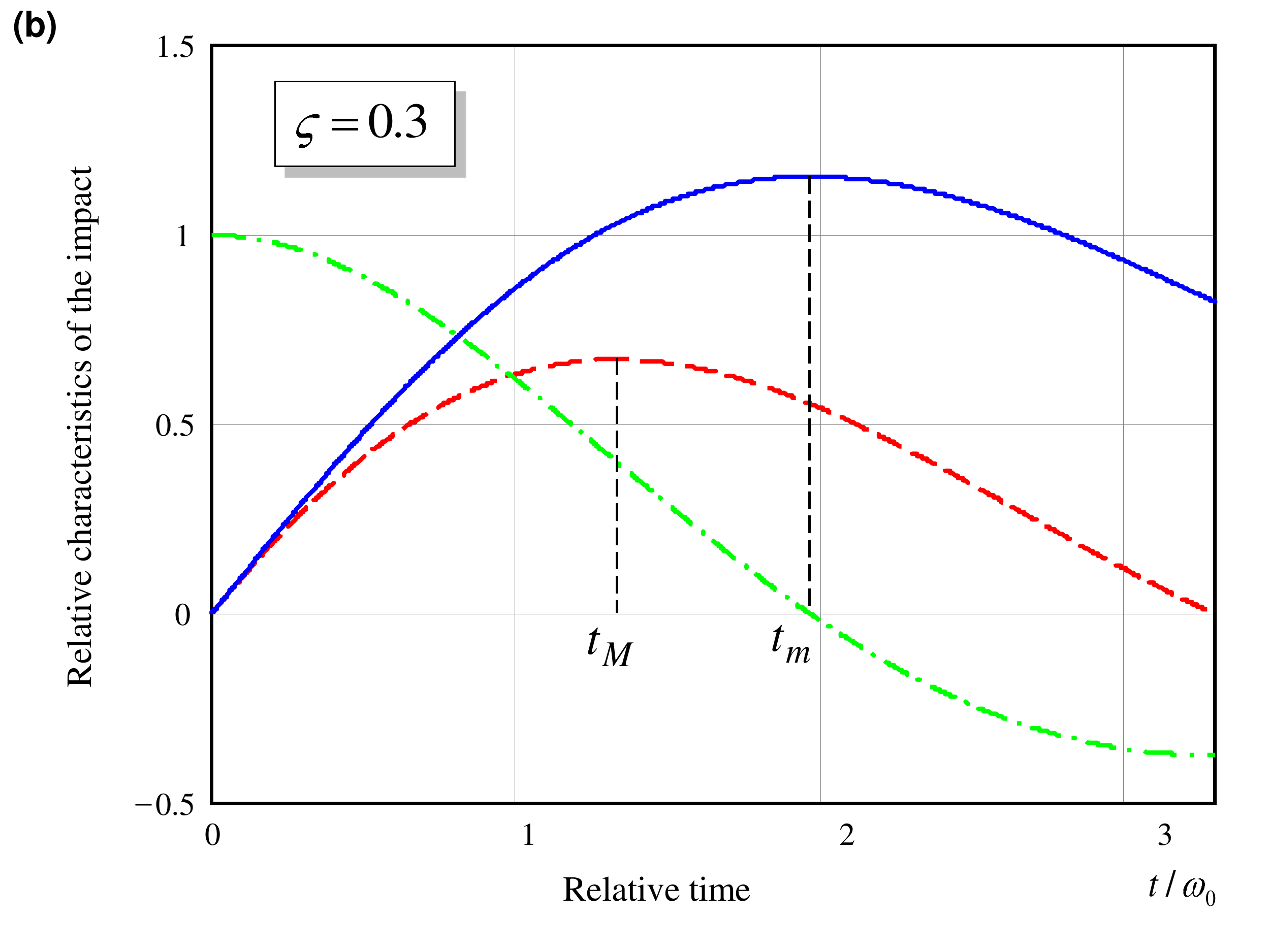}
    \includegraphics [scale=0.35]{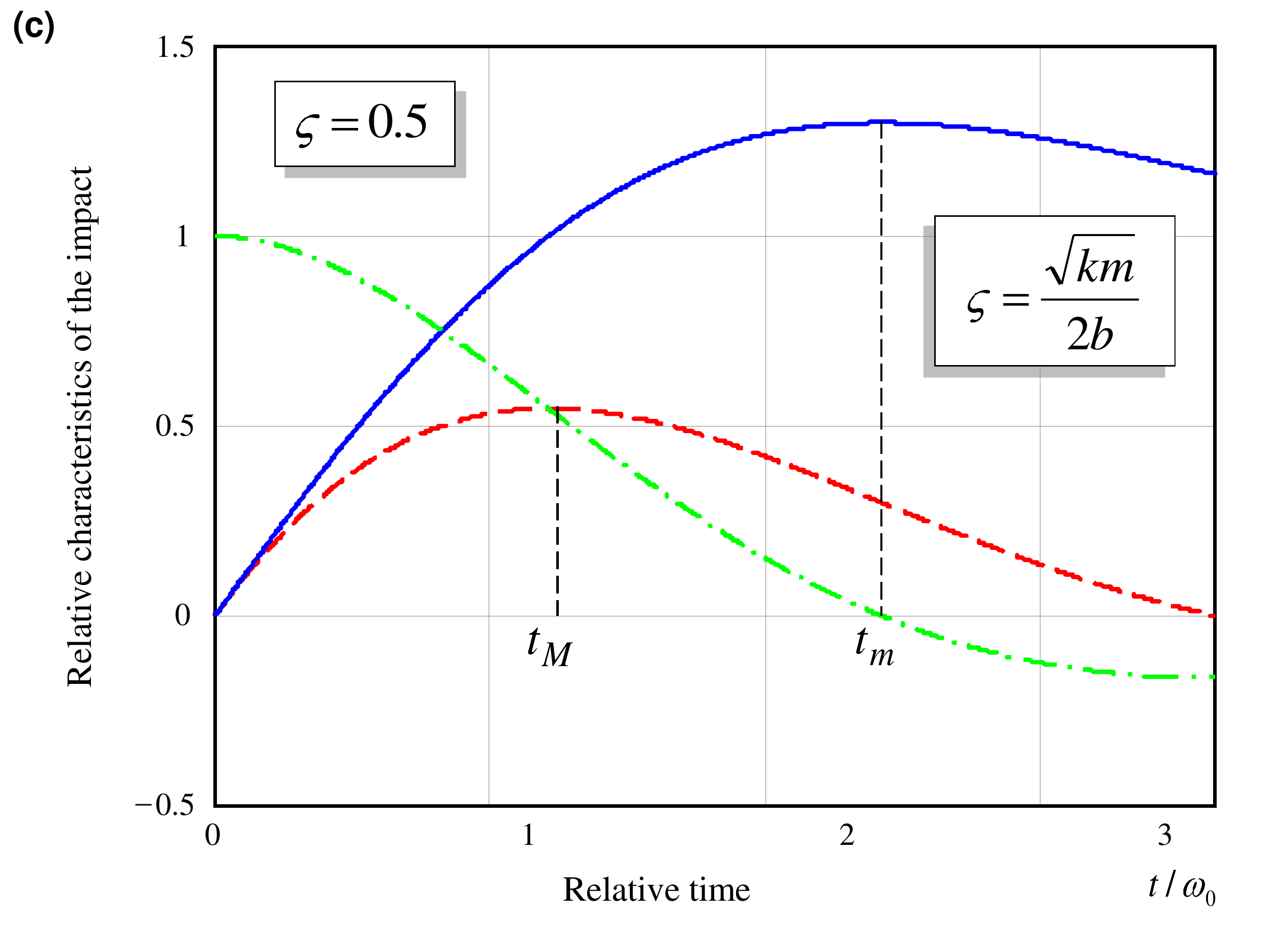}
    \caption{Viscoelastic Maxwell impact model. Behavior of the main impact variables with time for the following values of the damping ratio: $\zeta=0{.}1$ (a), $\zeta=0{.}3$ (b), $\zeta=0{.}5$ (c).}
    \label{Maxwell_exact.pdf}
\end{figure}

Fig.~\ref{Maxwell_exact.pdf} shows the behavior of the dimensionless quantities 
$\omega_0 x/v_0$, $F/(mv_0\omega_0)$, and $\dot{x}/v_0$ with respect to time. Observe that the time moment $t_m$, when the indenter displacement reaches its maximum, approaches the final moment of impact as the damping ratio $\zeta$ increases. 

According to Eq.~(\ref{1vI(2.5a)}), the peak value of the indenter penetration occurs at the instant 
\begin{equation}
t_m=\frac{\pi}{2\omega}\Bigl(1+\frac{2}{\pi}{\,\rm arcsin\,}\zeta\Bigr).
\label{1vI(2.10)}
\end{equation}
The substitution of the value (\ref{1vI(2.10)}) into Eqs.~(\ref{1vI(2.5)}) and (\ref{1vI(2.7)}) gives the maximum penetration
\begin{equation}
x_m=\frac{v_0}{\omega_0}\biggl(
2\zeta+\exp\biggl\{-\frac{\pi\zeta}{2\sqrt{1-\zeta^2}}\Bigl(1+\frac{2}{\pi}{\,\rm arcsin\,}\zeta\Bigr)
\biggr\}\biggr)
\label{1vI(2.11)}
\end{equation}
and the corresponding force 
\begin{equation}
F_m=\frac{kv_0}{\omega_0}
\exp\biggl\{-\frac{\pi\zeta}{2\sqrt{1-\zeta^2}}\Bigl(1+\frac{2}{\pi}{\,\rm arcsin\,}\zeta\Bigr)
\biggr\}
\label{1vI(2.11a)}
\end{equation}

From Eq.~(\ref{1vI(2.7)}), it follows that the peak value $F_M$ of the contact force occurs at the instant 
\begin{equation}
t_M=\frac{1}{\omega}{\,\rm atan\,}\frac{\sqrt{1-\zeta^2}}{\zeta}.
\label{1vI(2.12)}
\end{equation}

Substituting (\ref{1vI(2.12)}) into Eqs.~(\ref{1vI(2.7)}) and (\ref{1vI(2.5)}), we obtain the maximum contact force
\begin{equation}
F_M=\frac{kv_0}{\omega_0}
\exp\biggl\{-\frac{\zeta}{\sqrt{1-\zeta^2}}{\,\rm atan\,}\Bigl(
\frac{\sqrt{1-\zeta^2}}{\zeta}\Bigr)\biggr\}
\label{1vI(2.13)}
\end{equation}
and the corresponding displacement 
\begin{equation}
x_M=\frac{v_0}{\omega_0}\biggl(2\zeta+
(1-4\zeta^2)
\exp\biggl\{-\frac{\zeta}{\sqrt{1-\zeta^2}}{\,\rm atan\,}
\frac{\sqrt{1-\zeta^2}}{\zeta}\biggr\}\biggr).
\label{1vI(2.14)}
\end{equation}

\begin{figure}[h!]
    \centering
    \includegraphics [scale=0.35]{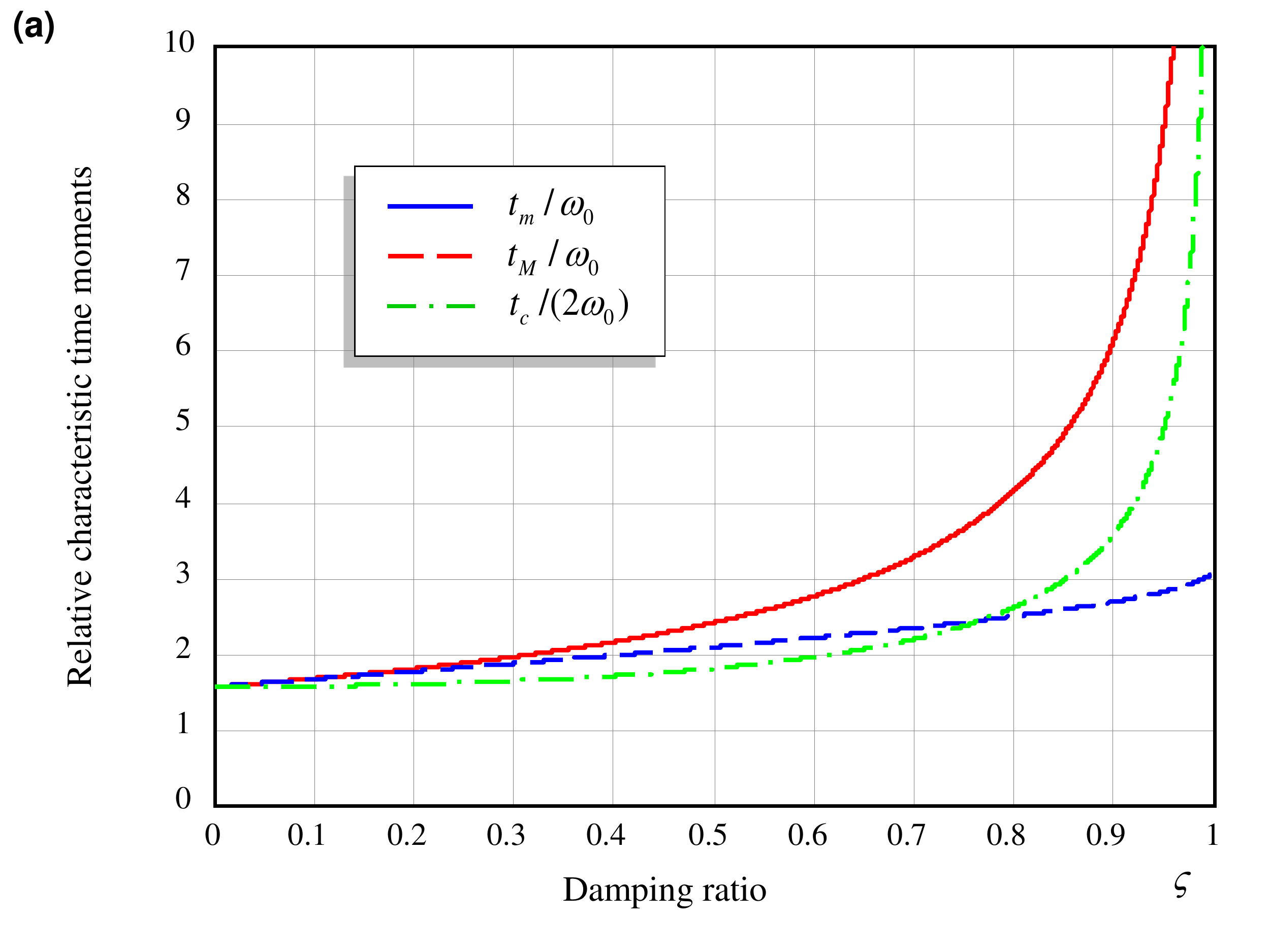}
    \includegraphics [scale=0.35]{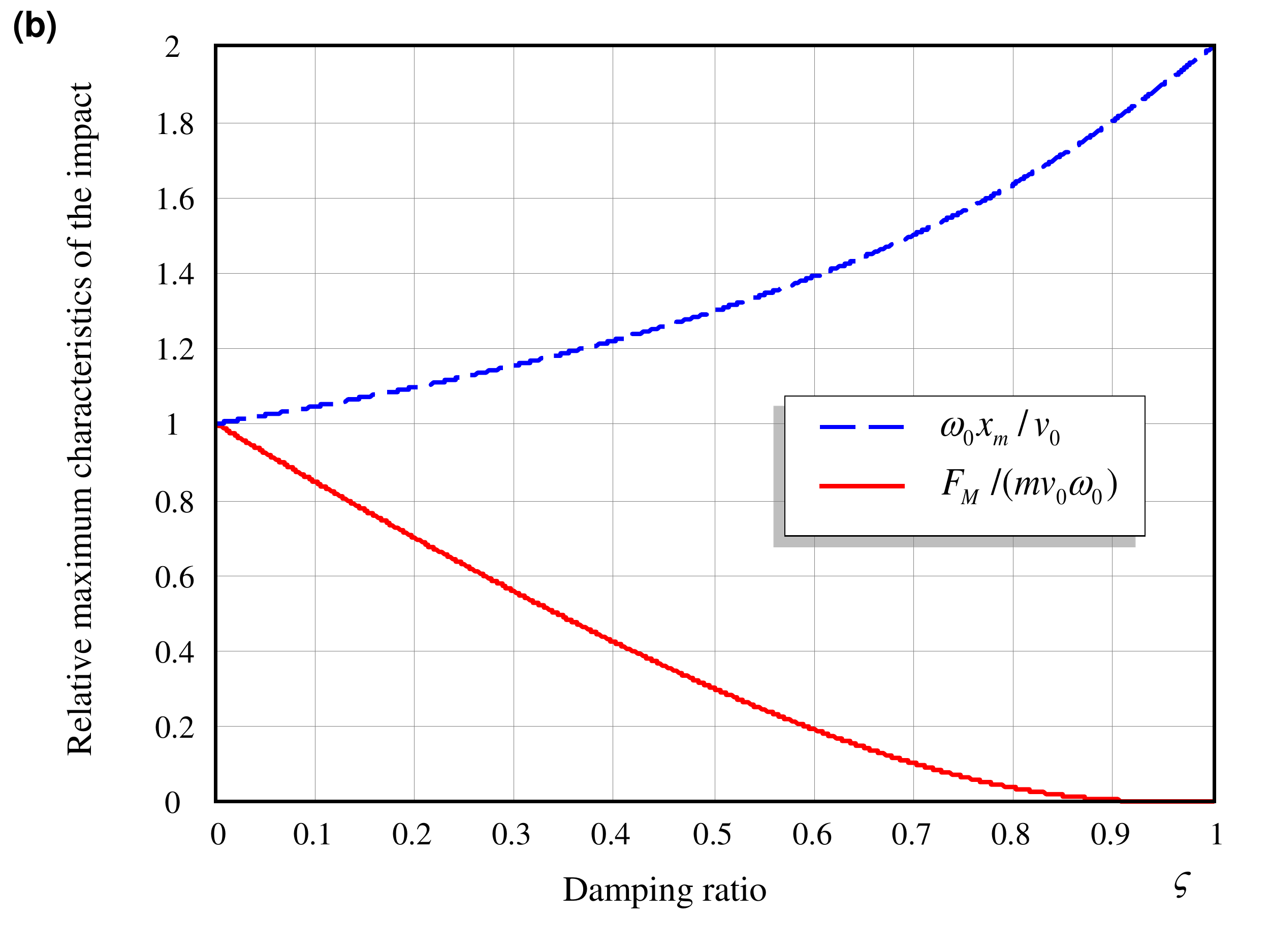}
    \caption{Viscoelastic Maxwell impact model. Behavior of the main impact parameters 
    $t_m$, $t_c$, $t_M$ (a) and $x_m$, $F_M$ (b) with the damping ratio.}
    \label{Maxwell_tm+tM.pdf}
\end{figure}

Fig.~\ref{Maxwell_tm+tM.pdf} shows the monotonic behavior of the dimensionless characteristic time moments $t_m/\omega_0$, $t_c/(2\omega_0)$, $t_M/\omega_0$ with the damping ratio $\zeta$. The variations of the relative maximum contact force $F_M/(mv_0\omega_0)$ and displacement $\omega_0 x_m/v_0$ are presented in
Fig.~\ref{Maxwell_tm+tM.pdf}b.

Finally, as it was observed \cite{ButcherSegalman2000}, although certain quantities of the Maxwell impact model are equivalent to the so-called half-period Kelvin\,--\,Voigt impact model, the inherent physics of these models are completely different. 

\section{Standard solid model}
\label{1dsSection3}

\begin{figure}[h!]
    \centering
    \includegraphics [scale=0.35]{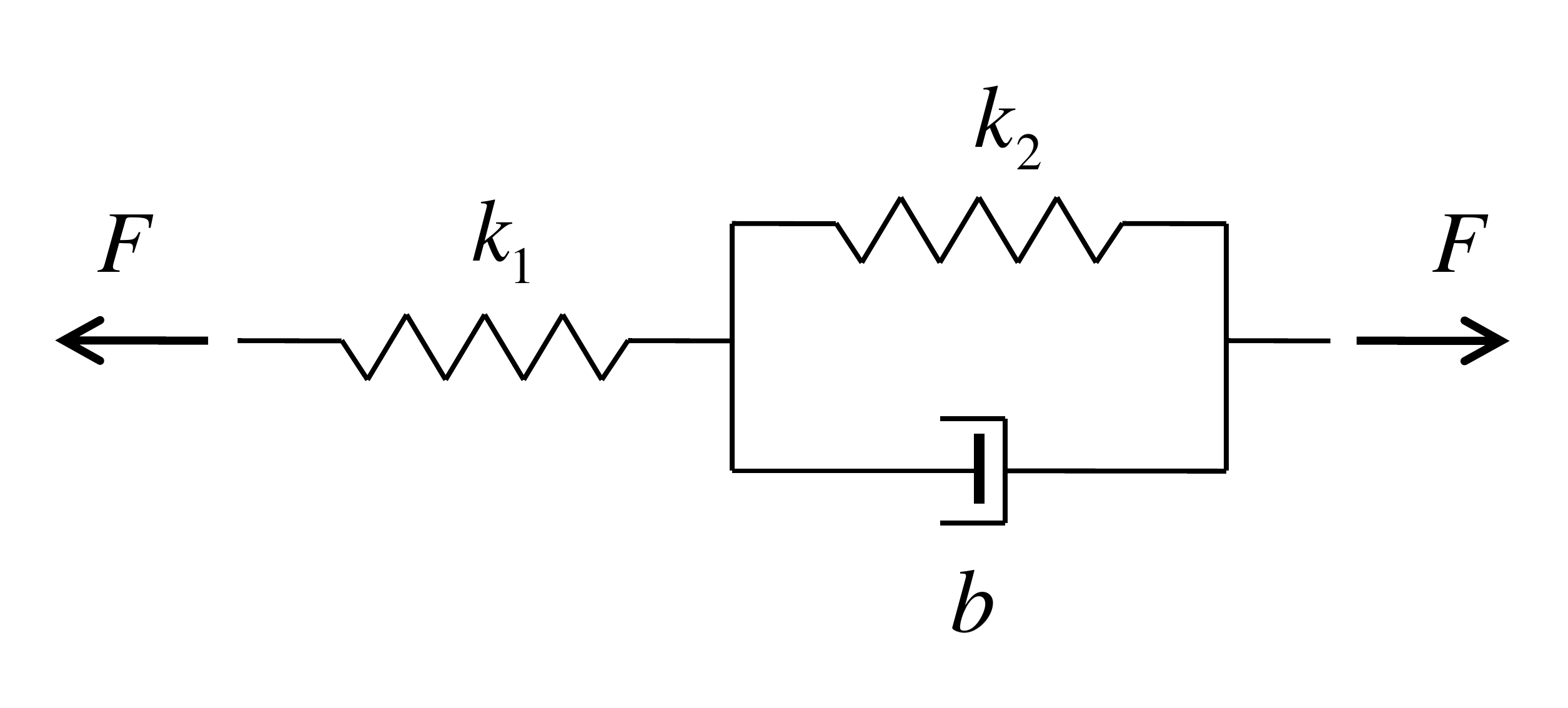}
    \caption{Standard solid model. Configuration based on the Kelvin--Voigt model.}
    \label{SSM_model-1.pdf}
\end{figure}

\begin{figure}[h!]
    \centering
    \includegraphics [scale=0.35]{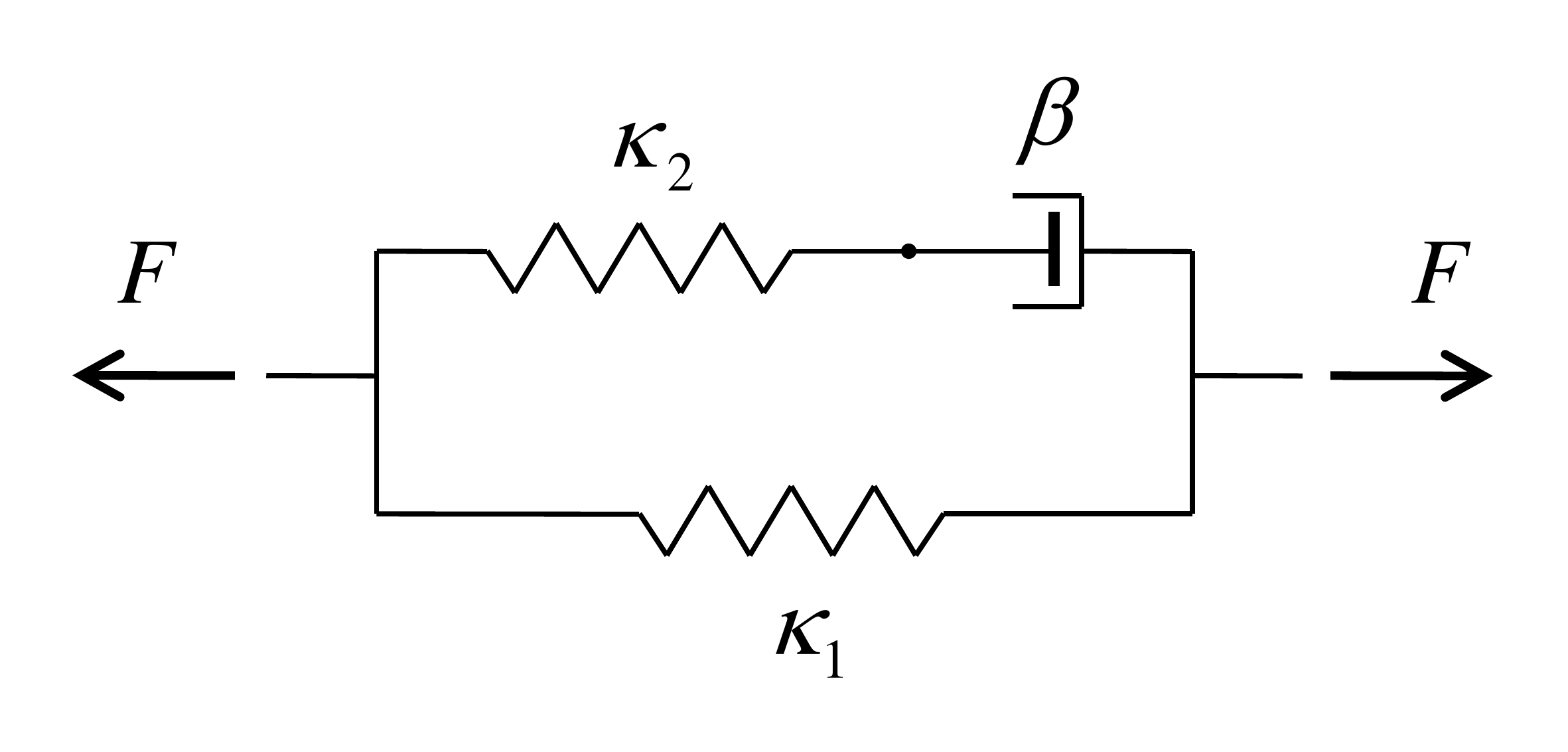}
    \caption{Standard solid model. Configuration based on the Maxwell model.}
    \label{SSM_model-2.pdf}
\end{figure}

There are two schematic representations of the standard linear solid model (Figs.~\ref{SSM_model-1.pdf} and \ref{SSM_model-2.pdf}). The force-displacement relationship is given by the following two equations:
\begin{equation}
(k_1+k_2)F+b\dot{F}=k_1 k_2 x+k_1 b\dot{x},
\label{1vI(3.1)}
\end{equation}
\begin{equation}
\varkappa_1 F+\beta\dot{F}=\varkappa_1 \varkappa_2 x+\beta(\varkappa_1+\varkappa_2)\dot{x}.
\label{1vI(3.2)}
\end{equation}

The instantaneous and long-term moduli are 
\begin{equation}
k_0=k_1=\varkappa_1+\varkappa_2,\quad
k_\infty=\varkappa_1=\frac{k_1 k_2}{k_1+k_2}.
\label{1vI(3.3)}
\end{equation}

The relaxation time is equal to 
\begin{equation}
\tau_R=\frac{b}{k_1+k_2}=\frac{\beta}{\varkappa_2}.
\label{1vI(3.4)}
\end{equation}

The differential equations (\ref{1vI(3.1)}) and (\ref{1vI(3.2)}) are equivalent to the force-displacement relationship 
\begin{equation}
F=\int\limits_0^t k(t-\tau)\frac{dx}{d\tau}(\tau)\,d\tau
\label{1vI(3.5)}
\end{equation}
with the relaxation stiffness
\begin{equation}
k(t)=k_\infty+(k_0-k_\infty)\exp\Bigl(-\frac{t}{\tau_R}\Bigr).
\label{1vI(3.6)}
\end{equation}

Let us also introduce the notation
\begin{equation}
\rho=\frac{k_\infty}{k_0}.
\label{1vI(3.7)}
\end{equation}
Note that $\rho\in(0,1)$ in view of (\ref{1vI(3.3)}).

The differential equation of impact
\begin{equation}
m\ddot{x}+F=0, \quad t\in[0,t_c],
\label{1vI(3.8)}
\end{equation}
where the contact force $F$ is determined by Eq.~(\ref{1vI(3.1)}), can be written as
\begin{equation}
\dddot{x}+\frac{(k_1+k_2)}{b}\ddot{x}+\frac{k_1}{m}\dot{x}+\frac{k_1 k_2}{mb}x=0.
\label{1vI(3.9)}
\end{equation}

By introducing the non-dimensional time 
\begin{equation}
\tau=\frac{t}{\tau_R},
\label{1vI(3.10)}
\end{equation}
Eq.~(\ref{1vI(3.9)}) can be reduced to the following equation:
\begin{equation}
x^{\prime\prime\prime}+x^{\prime\prime}+\Lambda x^\prime+\Lambda\rho x=0.
\label{1vI(3.11)}
\end{equation}
Here we introduced the notation 
\begin{equation}
\Lambda=\frac{k_0}{m}\tau_R^2.
\label{1vI(3.12)}
\end{equation}

The initial conditions for Eq.~(\ref{1vI(3.9)}) are as follows:
\begin{equation}
x(0)=0, \quad \dot{x}(0)=v_0,\quad \ddot{x}(0)=0.
\label{1vI(3.13)}
\end{equation}

The solution to the problem (\ref{1vI(3.9)}), (\ref{1vI(3.13)}) is given by the following formula:
\begin{eqnarray}
x(t) & = & \frac{\tau_R v_0}{\zeta_1[(\beta_1-\lambda_1)^2+\zeta_1^2]}
\biggl\{[(1-\beta_1)(\lambda_1-\beta_1)+\zeta_1^2]\sin\frac{\zeta_1 t}{\tau_R}
\nonumber \\
{} & { } & {}-\zeta_1(1-\lambda_1)\cos\frac{\zeta_1 t}{\tau_R}\biggr\}
\exp\Bigl(-\frac{\beta_1 t}{\tau_R}\Bigr)
+\frac{(1-\lambda_1)\tau_R v_0}{(\beta_1-\lambda_1)^2+\zeta_1^2}
\exp\Bigl(-\frac{\lambda_1 t}{\tau_R}\Bigr).
\label{1vI(3.14)}
\end{eqnarray}
Here, $-\lambda_1$ and $-(\beta_1\pm{\rm i}\zeta_1)$ are the roots of the characteristic equation corresponding to Eq.~(\ref{1vI(3.11)}). In other words, the following factorization takes place:
\begin{equation}
z^3+z^2+\Lambda z+\Lambda\rho =(z+\lambda_1)
(z^2+2\beta_1 z +\beta_1^2+\zeta_1^2).
\label{1vI(3.15)}
\end{equation}
The discriminant of the characteristic equation is 
\begin{equation}
D=4\Lambda(\Lambda^2+\rho)-\Lambda^2(1+18\rho-27\rho^2).
\label{1vI(3.16)}
\end{equation}

We underline that formula (\ref{1vI(3.14)}) holds true when $D>0$. In this case, we have
\begin{eqnarray}
\lambda_1 & = & \frac{1}{3}+\frac{C_1}{3}+\frac{(1-3\Lambda)}{3C_1},
\nonumber \\
\beta_1 & = & \frac{1}{3}-\frac{C_1}{6}-\frac{(1-3\Lambda)}{6C_1},
\nonumber \\
\zeta_1 & = & \frac{\sqrt{3}C_1}{6}-\frac{\sqrt{3}(1-3\Lambda)}{6C_1},
\nonumber
\end{eqnarray}
where
$$
C_1=\sqrt[3]{\frac{1}{2}(Q_1+2-9\Lambda+27\Lambda\rho)},\quad
Q_1=\sqrt{(2-9\Lambda+27\Lambda\rho)^2-4(1-3\Lambda)^3}.
$$

\section{Perturbation of the Kelvin--Voigt model}
\label{1dsSection4}

Taking into account (\ref{1vI(3.3)}) and (\ref{1vI(3.7)}), we rewrite Eq.~(\ref{1vI(3.1)}) in the following form:
\begin{equation}
k_\infty F+\rho(1-\rho)b\dot{F}=k_\infty^2 x+(1-\rho)k_\infty b\dot{x}.
\label{1vI(4.1)}
\end{equation}

Now, letting $\rho\to 0$, we arrive at the equation
\begin{equation}
F=k_\infty x+b\dot{x},
\label{1vI(4.2)}
\end{equation}
which coincides with Eq.~(\ref{1vI(1.0)}). Thus, for small values of $\rho$, the standard solid model (\ref{1vI(4.1)}) is a perturbation of the Kelvin--Voigt model (\ref{1vI(4.2)}).

Let us introduce the notation
\begin{equation}
\omega_0^2=\frac{k_\infty}{m},\quad 
\beta=\frac{b}{2m},\quad 
\eta=\frac{\beta}{\omega_0}.
\label{1vI(4.3)}
\end{equation}

Then, the parameters (\ref{1vI(3.4)}) and (\ref{1vI(3.12)}) can be evaluated as
\begin{equation}
\tau_R=2\eta\rho(1-\rho)\frac{1}{\omega_0},\quad 
\Lambda=4\eta^2\rho(1-\rho)^2.
\label{1vI(4.4)}
\end{equation}

In view of (\ref{1vI(4.4)}), the discriminant (\ref{1vI(3.16)}) and the roots of the characteristic equation (\ref{1vI(3.15)}) can be asymptotically evaluated as follows:
\begin{eqnarray}
D & = & 16\eta^2(1-\eta^2)\rho^2+O(\rho^3),\quad \rho\to 0,
\nonumber \\
\lambda_1 & = & 1-4\eta^2\rho+O(\rho^2),
\nonumber \\
\beta_1 & = & 2\eta^2\rho+O(\rho^2), \quad 
\zeta_1 = 2\eta\sqrt{1-\eta^2}\rho+O(\rho^2).
\nonumber
\end{eqnarray}

Consequently, we obtain the following asymptotic formulas for the impact duration, $t_c$, and the coefficient of restitution, $e_*$:
\begin{equation}
\omega_0 t_c \simeq \frac{2}{\sqrt{1-\eta^2}}{\,\rm atan\,}\frac{\sqrt{1-\eta^2}}{\eta}+
\rho\biggl\{4\eta-\frac{8\eta^2}{\sqrt{1-\eta^2}}{\,\rm atan\,}\frac{\sqrt{1-\eta^2}}{\eta}
\biggr\},
\label{1vI(4.5)}
\end{equation}
\begin{equation}
e_* \simeq \exp\biggl(-\frac{2\eta}{\sqrt{1-\eta^2}}{\,\rm atan\,}\frac{\sqrt{1-\eta^2}}{\eta}
\biggr)\biggl\{1+\frac{4\rho\eta}{\sqrt{1-\eta^2}}{\,\rm atan\,}\frac{\sqrt{1-\eta^2}}{\eta}
\biggr\}.
\label{1vI(4.6)}
\end{equation}

\begin{figure}[h!]
    \centering
    \includegraphics [scale=0.35]{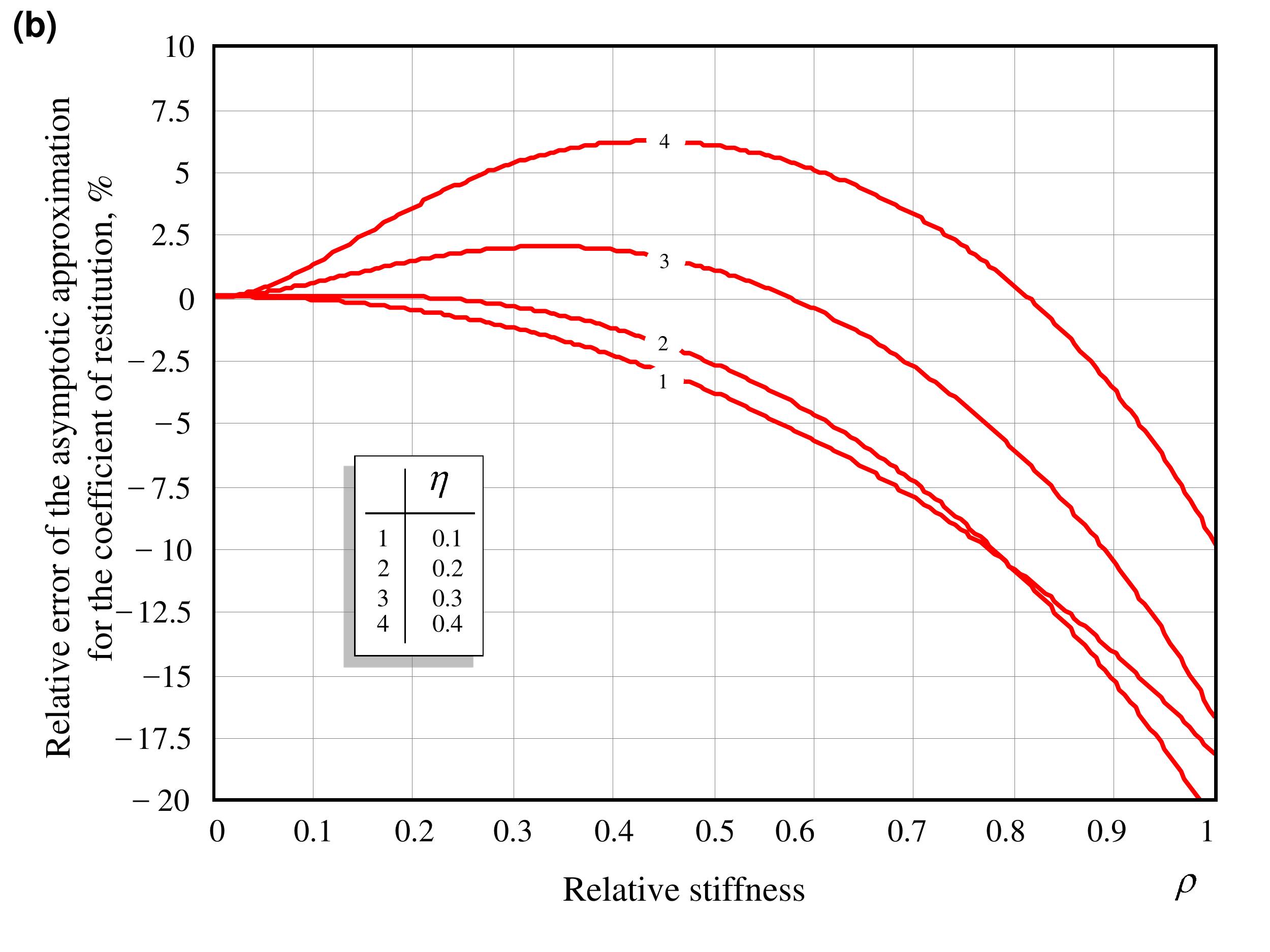}
    \includegraphics [scale=0.35]{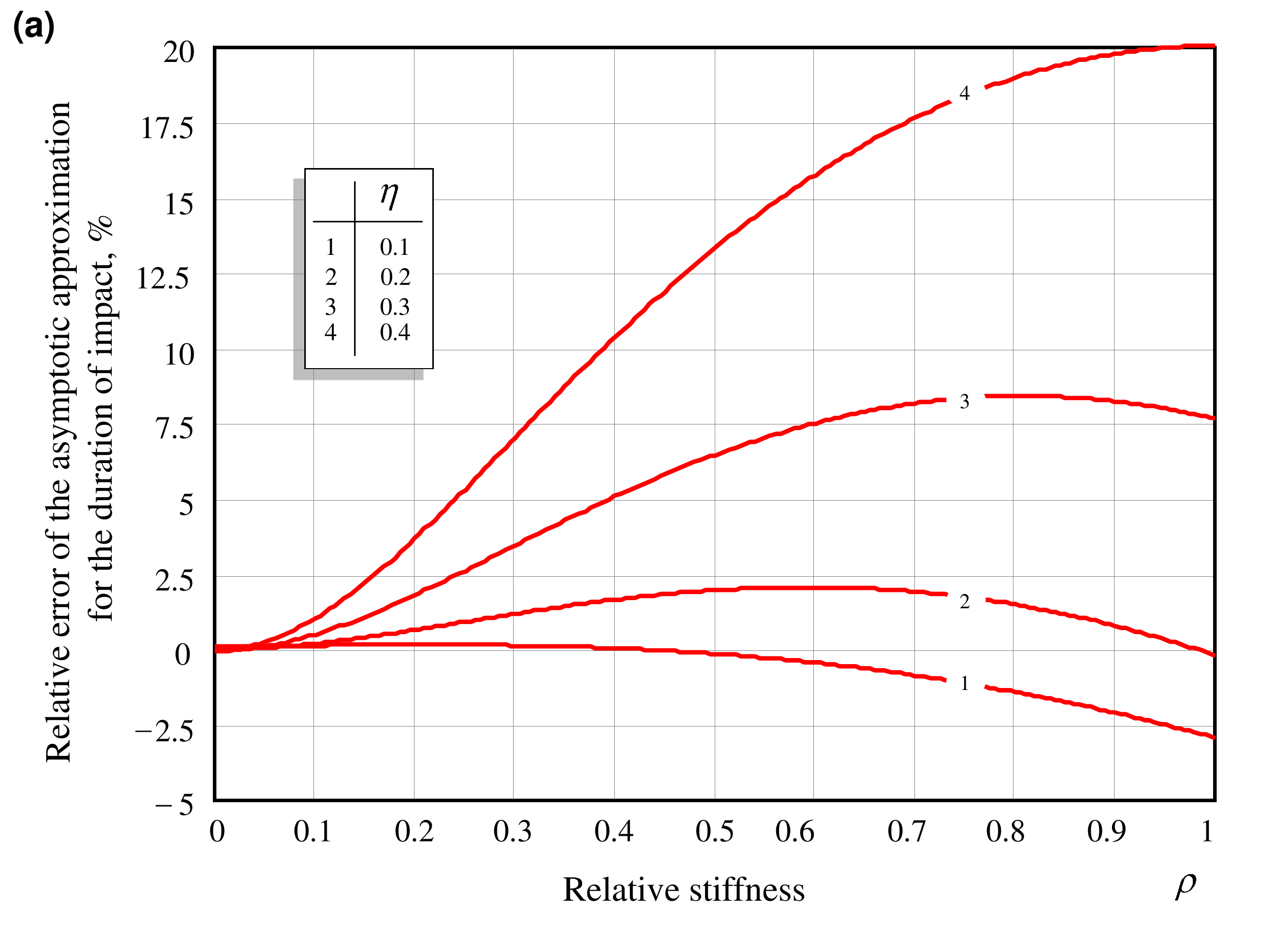}
        \caption{Perturbation of the Kelvin--Voigt model. Relative errors of the asymptotic formulas
                (\ref{1vI(4.5)}) and (\ref{1vI(4.6)}) for the duration of impact (a) and the coefficient of restitution (b).}
    \label{SSM_Kelvin-Voigt.pdf}
\end{figure}

The accuracy of the asymptotic approximations (\ref{1vI(4.5)}) and (\ref{1vI(4.6)}) is presented in 
Fig.~\ref{SSM_Kelvin-Voigt.pdf}. Note that the asymptotic formulas (\ref{1vI(4.5)}) and (\ref{1vI(4.6)}) are not uniformly valid as $\eta\to 1$.

\section{Perturbation of the Maxwell model}
\label{1dsSection5}

Now, taking into account (\ref{1vI(3.3)}) and (\ref{1vI(3.7)}), we rewrite Eq.~(\ref{1vI(3.1)}) as follows:
\begin{equation}
k_0 F+(1-\rho)b\dot{F}=\rho k_0^2 x+(1-\rho)k_0 b\dot{x}.
\label{1vI(5.1)}
\end{equation}

Again, by letting $\rho\to 0$, we obtain the limit equation
\begin{equation}
\frac{F}{b}+\frac{\dot{F}}{k_0}=\dot{x},
\label{1vI(5.2)}
\end{equation}
which coincides with Eq.~(\ref{1vI(2.1)}). Thus, for small values of $\rho$, the standard solid model (\ref{1vI(5.1)}) can be regarded as a perturbation of the Maxwell model (\ref{1vI(5.2)}).

Let us introduce the notation
\begin{equation}
\omega_0^2=\frac{k_0}{m},\quad 
\zeta=\frac{k_0}{2\omega_0 b}.
\label{1vI(5.3)}
\end{equation}

In view of (\ref{1vI(5.3)}), the parameters (\ref{1vI(3.4)}) and (\ref{1vI(3.12)}) can be evaluated as
\begin{equation}
\tau_R=\frac{1-\rho}{2\zeta\omega_0},\quad 
\Lambda=\frac{(1-\rho)^2}{4\zeta^2}.
\label{1vI(5.4)}
\end{equation}

Now, taking into account (\ref{1vI(5.4)}), we expand the discriminant (\ref{1vI(3.16)}) and the roots of the characteristic equation (\ref{1vI(3.15)}) as follows:
\begin{eqnarray}
D & = & \frac{(1-\zeta^2)}{16\zeta^6}-\frac{\rho(7\zeta^2-8\zeta^4+3)}{8\zeta^6}+O(\rho^2),\quad \rho\to 0,
\nonumber \\
\lambda_1 & = & \rho+O(\rho^2),\quad
\beta_1 = \frac{1}{2}-\frac{\rho}{2}+O(\rho^2),
\nonumber \\
\zeta_1 & = & \frac{\sqrt{1-\zeta^2}}{2\zeta}-
\rho\frac{\sqrt{1-\zeta^2}(2+\zeta^2)}{2\zeta(1-\zeta^2)}+O(\rho^2).
\nonumber
\end{eqnarray}

Consequently, we obtain the following asymptotic approximations for the impact duration, $t_c$, and the coefficient of restitution, $e_*$:
\begin{equation}
\omega_0 t_c \simeq \frac{\pi}{\sqrt{1-\zeta^2}}+
\frac{2\pi\rho\zeta^2}{(1-\zeta^2)^{3/2}},
\label{1vI(5.5)}
\end{equation}
\begin{equation}
e_* \simeq \exp\Bigl(-\frac{\pi\zeta}{\sqrt{1-\zeta^2}}\Bigr)
+4\rho\zeta^2
\biggl\{1+\biggl(1-\frac{\pi\zeta}{2(1-\zeta^2)^{3/2}}\biggr)
\exp\Bigl(-\frac{\pi\zeta}{\sqrt{1-\zeta^2}}\Bigr)\biggr\}.
\label{1vI(5.6)}
\end{equation}

\begin{figure}[h!]
    \centering
    \includegraphics [scale=0.35]{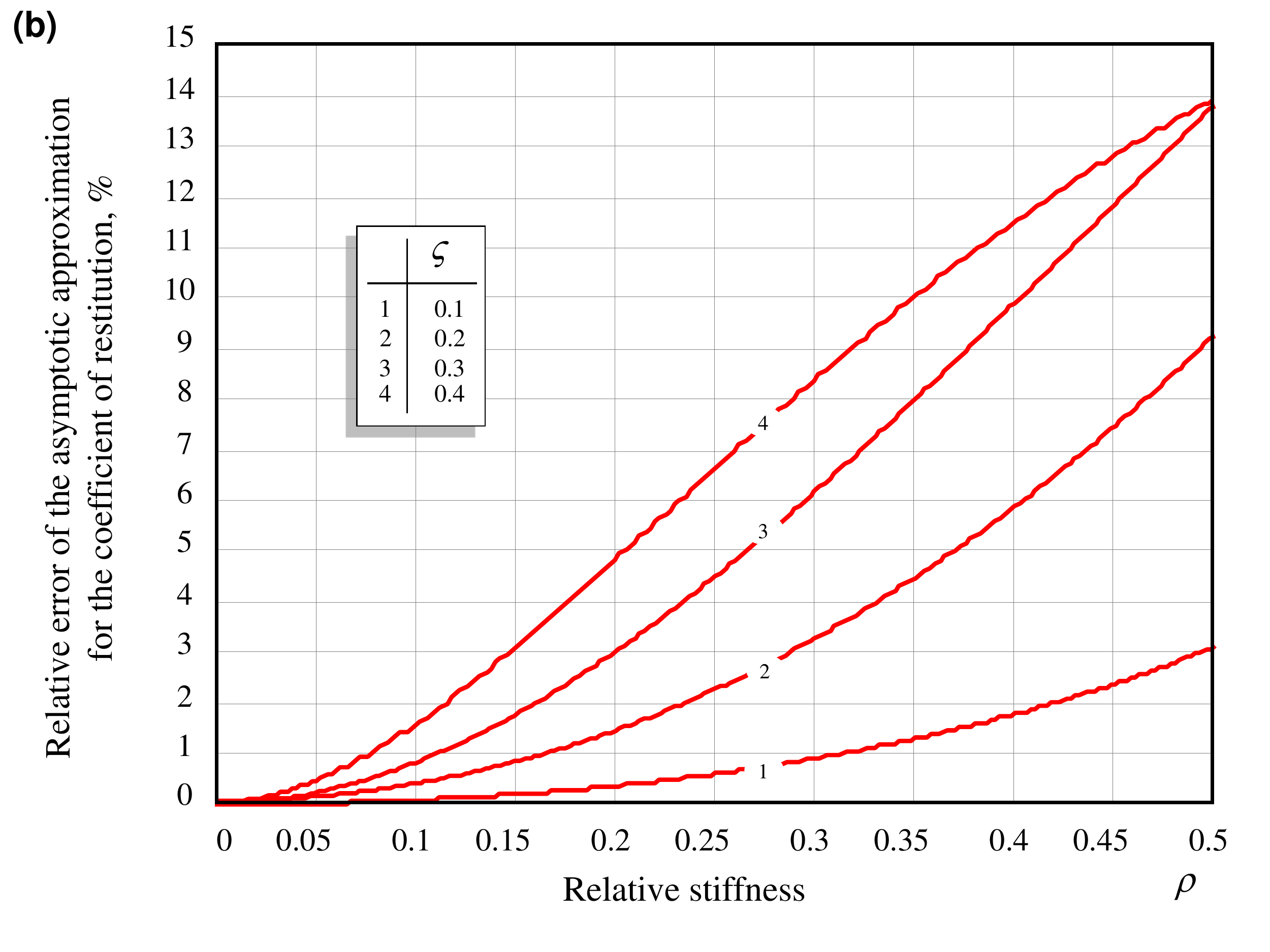}
    \includegraphics [scale=0.35]{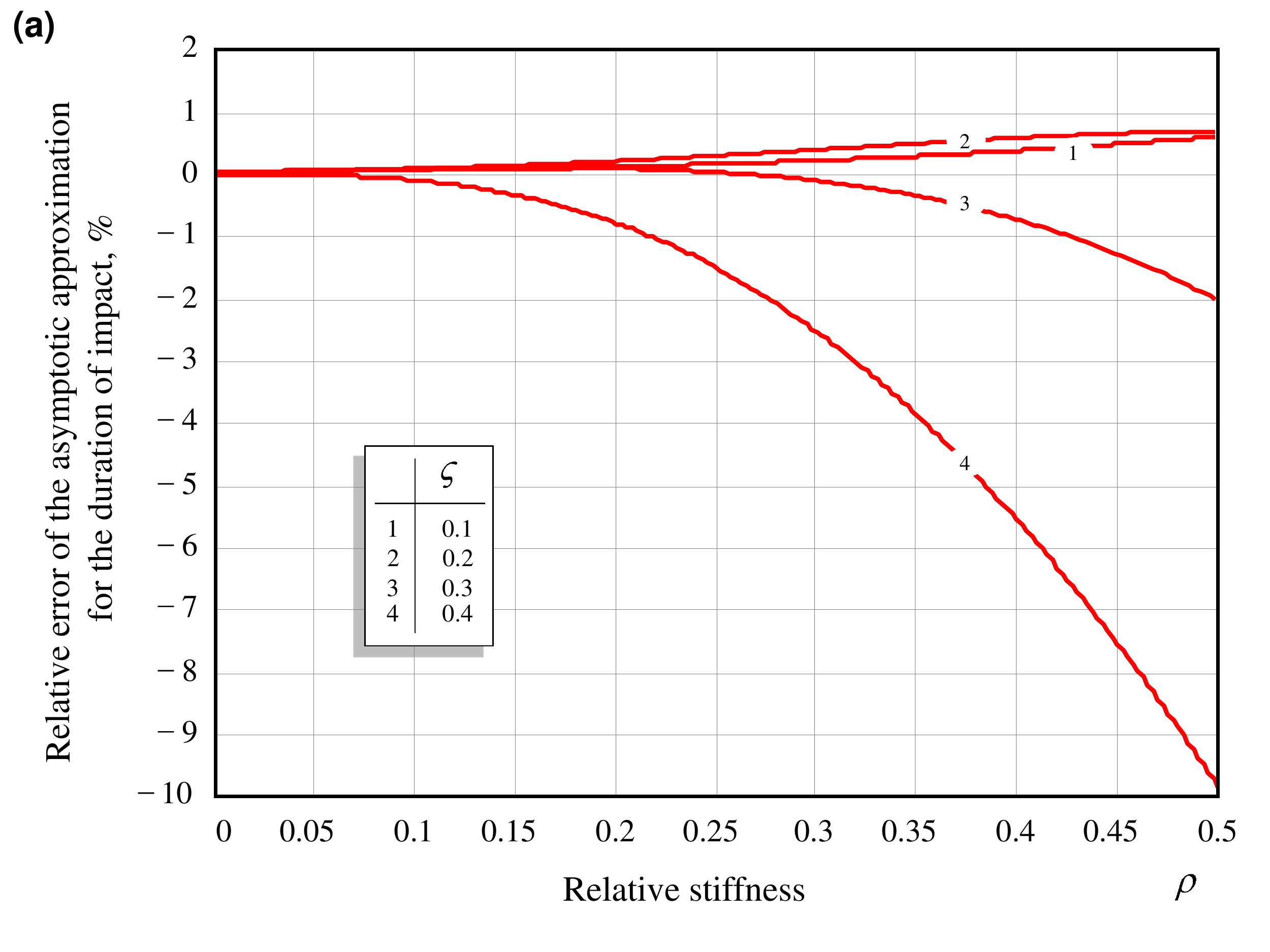}
        \caption{Perturbation of the Maxwell model. Relative errors of the asymptotic formulas
                (\ref{1vI(5.5)}) and (\ref{1vI(5.6)}) for the duration of impact (a) and the coefficient of restitution (b).}
    \label{SSM_Maxwell.pdf}
\end{figure}

The accuracy of the asymptotic approximations (\ref{1vI(5.5)}) and (\ref{1vI(5.6)}) is presented in 
Fig.~\ref{SSM_Maxwell.pdf}. Note that the asymptotic formulas (\ref{1vI(5.5)}) and (\ref{1vI(5.6)}) are not uniformly valid as $\eta\to 1$.

\section{Drop weight impact. Viscoelastic Kelvin--Voigt model}
\label{1dsSection05}

Due to Newton's second law, the differential equation of the drop weight impact has the form
\begin{equation}
m\ddot{x}+b\dot{x}+kx=mg, \quad t\in[0,t_c],
\label{1vI(05.1)}
\end{equation}
where $g$ is the gravitational acceleration. 

The initial conditions for Eq.~(\ref{1vI(05.1)}) are 
\begin{equation}
x(0)=0, \quad \dot{x}(0)=v_0.
\label{1vI(05.2)}
\end{equation}

The drop weight impact problem (\ref{1vI(05.1)}), (\ref{1vI(05.2)}) has the following solution:
\begin{equation}
x(t)=\frac{g}{\omega_0^2}\bigl(1-e^{-\beta t}\cos\omega t\bigr)+
\frac{1}{\omega}\Bigl(v_0-\frac{g\beta}{\omega_0^2}\Bigr)e^{-\beta t}\sin\omega t,
\label{1vI(05.3)}
\end{equation}
\begin{equation}
\dot{x}(t)=v_0 e^{-\beta t}\cos\omega t+
\frac{(g-\beta v_0)}{\omega}e^{-\beta t}\sin\omega t.
\label{1vI(05.4)}
\end{equation}
Here we used the notation (\ref{1vI(1.6)}).

According to Eqs.~(\ref{1vI(05.3)}), (\ref{1vI(05.4)}), the reaction force 
$F(x,\dot{x})=kx+b\dot{x}$ is given by
\begin{eqnarray}
\frac{F}{mv_0\omega_0} & = & \varepsilon_0\biggl\{1+e^{-\beta t}\biggl(
\frac{\eta}{\sqrt{1-\eta^2}}\sin\omega t-\cos\omega t\biggr)\biggr\}
\nonumber \\
{} & {} & {}+e^{-\beta t}\biggl(
\frac{1-2\eta^2}{\sqrt{1-\eta^2}}\sin\omega t+\frac{2\eta}{\sqrt{1-\eta^2}}\cos\omega t\biggr),
\label{1vI(05.5)}
\end{eqnarray}
where we introduced the notation
\begin{equation}
\varepsilon_0=\frac{g}{\omega_0 v_0}.
\label{1vI(05.5e)}
\end{equation}

The problem (\ref{1vI(05.1)}), (\ref{1vI(05.2)}) was studied in \cite{Ivanov1997}, where the existence of the parameter domain of ``plastic impact'' was established. This means that for any $\eta>0$, there exists a unique value of $\varepsilon_0^*$ such that for all $\varepsilon_0>\varepsilon_0^*$ we have $F(x,\dot{x})>0$ in the time interval $t\in(0,+\infty)$. The critical value $\varepsilon_0^*$ of the parameter $\varepsilon_0$ determines the critical value $v_0^*$ of the initial velocity $v_0$ below which there is no rebound effect. 

With the aim of application to the drop weight impact testing, we consider the problem (\ref{1vI(05.1)}), (\ref{1vI(05.2)}) for small values of the dimensionless parameter $\varepsilon_0$ and construct an asymptotic solution for the coefficient of restitution. 

Let $t_c^0$ and $e_*^0$ be the impact duration and the coefficient of restitution for the Kelvin--Voigt impact model, correspondingly. According to Eqs.~(\ref{1vI(1.9a)}) and (\ref{1vI(1.14)}), we have
\begin{equation}
t_c^0=\frac{2}{\omega_0\sqrt{1-\eta^2}}{\,\rm atan\,}\frac{\sqrt{1-\eta^2}}{\eta},\quad
e_*^0=\exp\biggl\{- \frac{2\eta}{\sqrt{1-\eta^2}}{\,\rm atan\,}\frac{\sqrt{1-\eta^2}}{\eta}\biggr\}.
\label{1vI(05.6)}
\end{equation}

Now, solving the transcendental equation $F(x,\dot{x})\bigr\vert_{t=t_c}=0$ by a perturbation method to terms of the first order inclusive, we obtain 
\begin{equation}
t_c\simeq t_c^0+\varepsilon_0\frac{(1+e_*^0)}{e_*^0\omega_0},
\label{1vI(05.7)}
\end{equation}
\begin{equation}
e_*\simeq e_*^0(1-2\varepsilon_0\eta),
\label{1vI(05.8)}
\end{equation}
where $t_c^0$ and $e_*^0$ are given by Eqs.~(\ref{1vI(05.6)}).

From the asymptotic formulas (\ref{1vI(05.7)}) and (\ref{1vI(05.8)}), it is clearly seen that the gravitational effect increases the duration of the impact process and decreases the coefficient of restitution. But it is more interesting to observe that the coefficient of restitution $e_*$ increases with velocity $v_0$, since the parameter $\varepsilon_0$ is inversely proportional to $v_0$ . That is why the effect of decrease in the coefficient of restitution in the drop weight impact test experimentally observed in \cite{Edelsten2010} for the velocity range $v_0\in(0{.}7,1{.}4)$~m/s and extrapolated for the low velocity region by means of the nonlinear Kelvin--Voigt model $F(x,\dot{x})=kx+c\vert x\vert\dot{x}$ with no account for the impactor weight cannot be explained by the linear viscoelastic Kelvin--Voigt model considered in this section.

\section{Drop weight impact. Viscoelastic Maxwell model}
\label{1dsSection06}

By applying the approach \cite{Stronge2000}, the differential equation of motion
$m\ddot{x}+F=mg$ with the initial conditions $x(0)=0$ and $\dot{x}(0)=v_0$ in view of the constitutive relationship (\ref{1vI(2.1)}) can be reduced to the following problem:
\begin{equation}
\dddot{x}+\frac{k}{b}\ddot{x}+\frac{k}{m}\dot{x}-\frac{kg}{b}=0,
\label{1vI(06.1)}
\end{equation}
\begin{equation}
x(0)=0, \quad \dot{x}(0)=v_0,\quad \ddot{x}(0)=g.
\label{1vI(06.2)}
\end{equation}

The drop weight impact problem (\ref{1vI(06.1)}), (\ref{1vI(06.2)}) has the following exact solution:
\begin{eqnarray}
x(t) & = & \frac{v_0}{\omega_0}e^{-\zeta\omega_0 t}\biggl\{
\frac{\omega_0}{\omega}(1-2\zeta^2-\varepsilon_0\zeta)\sin\omega t
-(2\zeta+\varepsilon_0)\cos\omega t\biggr\}
\nonumber \\
{} & {} & {}+\frac{v_0}{\omega_0}\bigl(2\zeta+\varepsilon_0(1-2\zeta\omega_0 t)\bigr),
\label{1vI(06.3)}
\end{eqnarray}
\begin{eqnarray}
\dot{x}(t) & = & v_0 e^{-\zeta\omega_0 t}\biggl\{
\frac{\omega_0}{\omega}(\zeta+\varepsilon_0)\sin\omega t+\cos\omega t\biggr\}
\nonumber \\
{} & {} & {}+2\zeta\varepsilon_0 v_0.
\label{1vI(06.4)}
\end{eqnarray}
Here we used the notation (\ref{1vI(2.6)}), (\ref{1vI(05.5e)}).

According to Eqs.~(\ref{1vI(06.3)}) and (\ref{1vI(06.4)}), the reaction force 
$F=mg-\ddot{x}$ is given by
\begin{equation}
\frac{F}{mv_0\omega_0}=\frac{\omega_0}{\omega} e^{-\zeta\omega_0 t}\biggl\{
(1+\varepsilon_0\zeta)\sin\omega t-\frac{\omega}{\omega_0}\varepsilon_0\cos\omega t\biggr\}
+\varepsilon_0
\label{1vI(06.5)}
\end{equation}

Now, let $t_c^0$ and $e_*^0$ be the impact duration and the coefficient of restitution for the Maxwell impact model, correspondingly. According to Eqs.~(\ref{1vI(2.8)}) and (\ref{1vI(2.9)}), we have
\begin{equation}
t_c^0=\frac{\pi}{\omega_0\sqrt{1-\zeta^2}},\quad
e_*^0=\exp\biggl(- \frac{\pi\zeta}{\sqrt{1-\zeta^2}}\biggr).
\label{1vI(06.6)}
\end{equation}

Applying a perturbation method, we find
\begin{equation}
t_c\simeq t_c^0+\varepsilon_0\frac{(1+e_*^0)}{e_*^0\omega_0},
\label{1vI(06.7)}
\end{equation}
\begin{equation}
e_*\simeq e_*^0-2\zeta\varepsilon_0,
\label{1vI(06.8)}
\end{equation}
where $t_c^0$ and $e_*^0$ are given by Eqs.~(\ref{1vI(06.6)}).

It is interesting to note that the asymptotic formulas (\ref{1vI(05.7)}) and (\ref{1vI(06.7)}) coincide. Clearly, the same conclusions can be drawn about the influence of the gravitational effect in the framework of the viscoelastic Maxwell drop weight impact model as those that were formulated in Section~\ref{1dsSection05}. 

\section{Short-time asymptotic solution of the indentation problem for a thin biphasic layer}
\label{1dsSection07}

We assume that the deformational behavior of articular cartilage is modeled in the framework of linear biphasic theory \cite{MowKueiLaiArmstrong1980}, which represents the biological tissue as a mixture consisting of a porous solid phase and a fluid phase (mobile interstitial water). The constitution equations for the solid and fluid phase stresses, $\mbox{$\boldsymbol\sigma$}^s$ and $\mbox{$\boldsymbol\sigma$}^f$, are given by
$$
\mbox{$\boldsymbol\sigma$}^s=-\phi^s p{\bf I}+\lambda^s {\rm tr}(\mbox{$\boldsymbol\varepsilon$}){\bf I}
+2\mu^s \mbox{$\boldsymbol\varepsilon$},\quad
\mbox{$\boldsymbol\sigma$}^f=-\phi^f p{\bf I}.
$$
Here, $\phi^f$ is the fluid volume fraction (porosity), $\phi^s=1-\phi^f$ is the solid volume fraction, $p$ is the true pressure of the fluid, $\lambda^s$ and $\mu^s$ are the Lam\'e constants, which together define the aggregate modulus $H_A=\lambda^s+2\mu^s$, $\mbox{$\boldsymbol\varepsilon$}$ is the strain tensor of the solid phase, and $\bf I$ is the identity tensor. Note that the fluid phase is assumed to be intrinsically incompressible and inviscid. 

The continuity equation for the mixture and the momentum equations for each phase are given by 
$$
{\rm div}(\phi^s{\bf v}^s+\phi^f{\bf v}^f)=0,
$$
$$
{\rm div}\mbox{$\boldsymbol\sigma$}^s+\frac{(\phi^f)^2}{\kappa}({\bf v}^f-{\bf v}^s)={\bf 0}, \quad
{\rm div}\mbox{$\boldsymbol\sigma$}^f-\frac{(\phi^f)^2}{\kappa}({\bf v}^f-{\bf v}^s)={\bf 0},
$$
where ${\bf v}^s$ and ${\bf v}^f$ are the solid and fluid velocities, respectively, and $\kappa$ is the permeability of the solid phase. 

Let us consider an axisymmetric contact problem for a thin biphasic layer indented without friction by a rigid impermeable cylindrical indenter. It is assumed that the contact radius $a$ is much larger than the cartilage layer thickness $h$ (i.e., $h/a\ll 1$). In this case, according to \cite{Ateshian1994}, the vertical displacements $w(r,t)$ of the boundary points of the articular cartilage tissue at the contact zone can be approximated by the following asymptotic formula:
\begin{equation}
w(r,t) = \frac{h^3}{3\mu_s}\biggl\{
\frac{1}{3r}\frac{\partial}{\partial r}\Bigl(r\frac{\partial P}{\partial r}(r,t)\Bigr)
+\frac{\mu_s \kappa}{h^2}\int\limits_0^t
\frac{1}{r}\frac{\partial}{\partial r}\Bigl(r\frac{\partial P}{\partial r}(r,\tau)\Bigr)d\tau\biggr\}.
\label{1vI(9.1)}
\end{equation}
Here, $P(r,t)$ is the contact pressure. It is assumed that the cartilage layer is bonded to a rigid impermeable substrate, that is there is no solid displacement at the cartilage-bone interface and no fluid flow through the bone \cite{Ateshian1994}. 

In view of (\ref{1vI(9.1)}), the contact condition that the boundary points of the cartilage layer acquire a constant vertical displacement $-\delta_0(t)$ (due to the action of the indenter) can be written as 
\begin{equation}
w(r,t) = -\delta_0(t), \quad r\leq a.
\label{1vI(9.2)}
\end{equation}

The substitution of (\ref{1vI(9.1)}) into Eq.~(\ref{1vI(9.2)}) results in an integro-differential equation 
\begin{equation}
\frac{1}{r}\frac{\partial}{\partial r}\Bigl(r\frac{\partial P}{\partial r}(r,t)\Bigr)
+\frac{3\mu_s \kappa}{h^2}\int\limits_0^t
\frac{1}{r}\frac{\partial}{\partial r}\Bigl(r\frac{\partial P}{\partial r}(r,\tau)\Bigr)d\tau
=-\frac{3\mu_s}{h^3}\delta_0(t),
\label{1vI(9.2a)}
\end{equation}
which requires imposing a suitable boundary condition at the edge of the contact zone, i.e., at $r=a$ .

In order to impose the mentioned boundary condition, we note that at the initial moment of contact $t=0$, formula (\ref{1vI(9.1)}) simplifies as follows:
\begin{equation}
w(r,0) = \frac{h^3}{3\mu_s}\frac{1}{3r}\frac{\partial}{\partial r}\Bigl(r\frac{\partial P}{\partial r}(r,0)\Bigr).
\label{1vI(9.3)}
\end{equation}
Comparing formula (\ref{1vI(9.3)}) with the known asymptotic solutions for thin elastic layers \cite{Barber1990,Chadwick2002,ArgatovMishuris2011ve}, we conclude that the instantaneous deformational response of a thin biphasic layer coincides with the response of a thin bonded incompressible elastic layer.
Thus, by this analogy, we will require that $P(r,t)\to 0$ as $r\to a$, that is the contact pressure is assumed to vanish at the edge of the contact area. 

As a result of integration of Eq.~(\ref{1vI(9.2a)}) with respect to the radial coordinate, we arrive at the following integral equation:
\begin{equation}
P(r,t)+\frac{3\mu_s \kappa}{h^2}\int\limits_0^t P(r,\tau)\,d\tau
=\frac{3\mu_s}{4h^3}\delta_0(t)(a^2-r^2).
\label{1vI(9.4)}
\end{equation}

Now, in order to derive the relationship between the indenter displacement and the contact force
$$
F(t)=2\pi\int\limits_0^a P(r,t)r\,dr,
$$
we multiply both sides of Eq.~(\ref{1vI(9.4)}) by $2\pi r$ and after that we integrate the equation obtained with respect to $r$ from $0$ to $a$. As a results of this operation, we get
\begin{equation}
F(t)+\frac{3\mu_s \kappa}{h^2}\int\limits_0^t F(\tau)\,d\tau
=\frac{3\mu_s a^4}{16h^3}\delta_0(t).
\label{1vI(9.5)}
\end{equation}

Further, by inverting the Volterra integral operator on the right-hand side of Eq.~(\ref{1vI(9.5)}), we obtain
\begin{equation}
F(t)=\frac{3\mu_s a^4}{16h^3}\biggl\{\delta_0(t)-
\chi\int\limits_0^t e^{-\chi(t-\tau)}\delta_0(\tau)\,d\tau\biggr\},
\label{1vI(9.6)}
\end{equation}
where we introduced the shorthand notation
\begin{equation}
\chi=\frac{3\mu_s \kappa}{h^2}.
\label{1vI(9.6b)}
\end{equation}

Finally, after integrating by parts, Eq.~(\ref{1vI(9.6)}) yields 
\begin{equation}
F(t)=\frac{3\mu_s a^4}{16h^3}\biggl\{\delta_0(0)+
\int\limits_0^t e^{-\chi(t-\tau)}\frac{d\delta_0}{d\tau}(\tau)\,d\tau\biggr\}.
\label{1vI(9.7)}
\end{equation}

In impact problems, under the assumption that 
$$
\delta_0(0)=0,
$$
the force-displacement relationship (\ref{1vI(9.7)}) takes the form
\begin{equation}
F(t)=\frac{3\mu_s a^4}{16h^3}\int\limits_0^t 
\exp\Bigl\{-\frac{(t-\tau)}{\tau_R}\Bigr\}\frac{d\delta_0}{d\tau}(\tau)\,d\tau\biggr\}.
\label{1vI(9.8)}
\end{equation}
Here we introduced the notation $\tau_R=1/\chi$. In view of (\ref{1vI(9.6b)}), we have 
\begin{equation}
\tau_R=\frac{h^2}{3\mu_s \kappa},
\label{1vI(9.9)}
\end{equation}
while comparing (\ref{1vI(9.9)}) with (\ref{1vI(2.2)}), we get the stiffness coefficient
\begin{equation}
k=\frac{3\mu_s a^4}{16h^3}.
\label{1vI(9.90)}
\end{equation}

It should be emphasized that Eq.~(\ref{1vI(9.8)}) represents a short-time asymptotic approximation, which is valid for moments of time $t$ such that $H_A\kappa t/h^2\ll 1$. For typical human cartilage material properties, $H_A=0{.}5$~MPa and $\kappa=2\times 10^{-15}$ ${\rm m}^4/{\rm Ns}$. Thus, assuming a typical cartilage thickness $h=1$~mm, we get $h^2/(H_A\kappa)=10^3$~s; thus, the asymptotic model (\ref{1vI(9.8)}) certainly remains valid for up to 100~s, which is well in the range of usual values of impact durations. 

Comparing Eq.~(\ref{1vI(9.8)}) with Eq.~(\ref{1vI(2.2)}), we see that the short-time deformational response of a thin biphasic layer bonded to a rigid impermeable substrate under the action of a frictionless flat-ended indenter is mathematically equivalent to that of a thin incompressible layer following the Maxwell viscoelastic model. Note that the Maxwell's model based perturbation model considered in Section~\ref{1dsSection5} could be useful for modeling the impact response of articular cartilage (or artificial tissues for its replacement) in the whole time range, i.e. in the short-, medium- and long-time range.

We also emphasize that the biphasic model is not equivalent to a viscoelastic model, because the biomechanical response of a poroelastic material such as articular cartilage is crucially dependent on the boundary conditions for the sample. In particular, viscoelastic equivalents of the deformational response of an articular cartilage sample subjected to the same simple loading protocols in confined and unconfined conditions will be essentially different, especially in the short-time range. Thus, in comparing experimental results from different sources, a particular attention should be paid to the fixation conditions for tissue samples. 

\section{Key features of non-linear impact}
\label{1dsSection10}

To illustrate the application of the developed linear theory of viscoelastic impact, let us analyze the experimental data obtained in \cite{BurginAspden2008}
for drop-weight impact testing (with the impactor mass $m=100$~g) of isolated bovine articular cartilage samples of 5~mm diameter (correspondingly, the radius of the samples is $a=2{.}5$~mm). 
In \cite{BurginAspden2008}, the force data, $F(t)$, were converted to engineering stress, $\sigma(t)$, by dividing them by the original cross-section area of the sample, $\pi a^2$, i.\,e.,
$$
\sigma(t)=\frac{F(t)}{\pi a^2}.
$$
The effective strain, $\epsilon(t)$, was evaluated by dividing the measured impactor displacement, $x(t)$, by the sample thickness, $h$, which is assumed to be $0{.}5\pm 0{.}11$~mm, as follows:
$$
\epsilon(t)=\frac{x(t)}{h}.
$$
(Here, stress and strain are assumed to be positive in compression.) 
The stress-strain relationship was differentiated to obtain the incremental dynamic modulus
$$
E_{\rm dyn}=\frac{d\sigma}{d\epsilon}.
$$
The maximum incremental dynamic modulus, $E_{\rm max}$, was found, and the modulus $E_{10}$ at stresses of 10~MPa was determined to enable comparison of dynamic moduli at constant value of stress. The initial impact velocity was calculated from the drop height, $h_0$, by the well-known formula
$v_0=\sqrt{2gh_0}$.

The incremental dynamic modulus can be evaluated as a function of time in the form
\begin{equation}
E_{\rm dyn}(t)=\frac{\dot{\sigma}(t)}{\dot{\epsilon}(t)}=\frac{h}{\pi a^2}\frac{\dot{F}(t)}{\dot{x}(t)}.
\label{1vI(E.1)}
\end{equation}
In the case of the Maxwell model (see, Section~\ref{1dsSection2}), we will have
\begin{equation}
\frac{\dot{F}(t)}{\dot{x}(t)}=k\frac{\cos\omega t-(\zeta\omega_0/\omega)\sin\omega t
}{\cos\omega t+(\zeta\omega_0/\omega)\sin\omega t},
\label{1vI(E.2)}
\end{equation}
where $k$ is the stiffness coefficient. 

First of all, observe that in view of (\ref{1vI(E.1)}) and (\ref{1vI(E.2)}), the variation of $E_{\rm dyn}(t)$ does not depend on the impact velocity $v_0$. In other words, the time variation of the incremental dynamic stiffness in the linear viscoelastic impact tests remains the same for different initial impact velocities. We emphasize that this conclusion is valid for a general linear viscoelastic law. Second, from (\ref{1vI(E.1)}) and (\ref{1vI(E.2)}), it follows that the value of $E_{\rm dyn}(t)$ gradually decreases to zero with increasing contact force $F(t)$ (when $\dot{F}(t)>0$). Thus, we arrive at the formula 
\begin{equation}
E_{\rm max}=E_{\rm dyn}(0). 
\label{1vI(E.0)}
\end{equation}

Further, in order to evaluate $E_{10}$, we need first solve the equation 
\begin{equation}
F(t_{10})=\pi a^2\sigma_{10},
\label{1vI(E.3)}
\end{equation}
where $\sigma_{10}=10$~MPa. In view of (\ref{1vI(2.7)}), Eq.~(\ref{1vI(E.3)}) takes the form
\begin{equation}
\exp(-\zeta\omega_0 t_{10})\sin\omega t_{10}=\frac{\pi a^2\omega\sigma_{10}}{kv_0}.
\label{1vI(E.4)}
\end{equation}
Here, $\zeta$, $\omega_0$, and $\omega$ are independent of $v_0$, and are determined by formulas 
(\ref{1vI(2.6)}).

Now, from (\ref{1vI(E.4)}), it is seen that the value of the time moment $t_{10}$ depends on the initial velocity $v_0$. Thus, the Maxwell impact model (and generally speaking, any linear viscoelastic model of impact) predicts that the value of $E_{10}$ increases with increasing impact velocity~$v_0$.

\begin{table}[h!]
\caption{Impact parameters for isolated bovine articular cartilage samples \cite{BurginAspden2008}}
\label{1vITable1}       
\smallskip
\begin{tabular}{cccccccc}
\hline
$h_0$ (mm) & $v_0$ (m/s) & $E_{\rm max}$ (MPa) & $E_{10}$ (MPa) & $\sigma_{\rm max}$ (MPa) & $\epsilon_{\max}$ & $e_*$ & $\Delta m$ (\%)  \\
\hline\noalign{\smallskip}
25 & $0{.}70$ & $86\pm 22$ & $75\pm 13$ & $15{.}6\pm 2{.}9$ & $0{.}48\pm 0{.}06$ & $0{.}64\pm 0{.}08$ & $2{.}2$ \\
50 & $0{.}99$ & $100\pm 32$ & $71\pm 16$ & $24{.}5\pm 3{.}5$ & $0{.}60\pm 0{.}13$ & $0{.}46\pm 0{.}14$ & $2{.}5$ \\
80 & $1{.}25$ & $118\pm 33$ & $73\pm 12$ & $34{.}2\pm 5{.}0$ & $0{.}62\pm 0{.}11$ & $0{.}47\pm 0{.}05$ & $5{.}7$ \\
100 & $1{.}40$ & $128\pm 28$ & $72\pm 13$ & $40{.}5\pm 4{.}6$ & $0{.}68\pm 0{.}09$ & $0{.}41\pm 0{.}08$ & $9{.}9$ \\
\hline
\end{tabular}
\end{table}

Table~\ref{1vITable1} shows that the impact testing \cite{BurginAspden2008} was performed in the non-linear regime with maximum compressive strains of 50--60\%. That is why, the prediction of the linear impact model concerning $E_{\rm max}$ are not fulfilled. Furthermore, the linear theories of impact predict that the maximum contact force $F_M$ (correspondingly, the maximum contact stress $\sigma_{\rm max}=F_M/(\pi a^2)$) and the maximum displacement $x_m$ (correspondingly, the maximum strain $\epsilon_{\rm max}=x_m/h$) are proportional to $v_0$. On the other hand, the data from Table~1 show that the ratio $\sigma_{\rm max}/\epsilon_{\rm max}$ increases with increasing $v_0$. This fact also clearly indicates the non-linear deformational behavior of cartilage at high level of strain. Note here that the ratio $\sigma_{\rm max}/\epsilon_{\rm max}$ is ralted to the so-called pulsatile dynamic modulus (see, in particular, \cite{Argatov2012sine}).

Concerning the coefficient of restitution $e_*$ note that it is not constant, as it would be if the cartilage deformation were described by the Maxwell model (see formula (\ref{1vI(2.9)})). 

Finally, the last column of Table~\ref{1vITable1} gives the values of percentage increase in mass of each sample after 24~h immersed in PBS following impact loading. This is indicative of increasing amounts of damage in the cartilage samples \cite{BurginAspden2008}.

\section{Discussion}
\label{1dsSectionD}

Consider now the general case of linear viscoelastic force-displacement relationship
\begin{equation}
F=\int\limits_0^t k(t-s)\frac{dx}{ds}(s)\,ds
\label{1vI(7.1)}
\end{equation}
with the relaxation stiffness
$$
k(t)=k_0\Psi\Bigl(\frac{t}{\tau_R}\Bigr).
$$
Here, $k_0$ is the initial stiffness, $\tau_R$ is the characteristic relaxation time, $\Psi(\tau)$ is the dimensionless relaxation function with $\tau$ being a dimensionless independent time-like variable.

Making use of the change of variables
\begin{equation}
t=\tau_R\tau,\quad x=v_0\tau_R\xi,
\label{1vI(7.2)}
\end{equation}
we transform the impact equation $m\ddot{x}+F=0$ and the initial conditions $x(0)=0$, $\dot{x}=v_0$ into the following problem:
\begin{equation}
\xi^{\prime\prime}+\alpha\int\limits_0^\tau \Psi(\tau-\sigma)\frac{d\xi}{d\sigma}(\sigma)\,d\sigma=0,
\label{1vI(7.3)}
\end{equation}
\begin{equation}
\xi(0)=0,\quad \xi^\prime(0)=1.
\label{1vI(7.4)}
\end{equation}
Here prime denotes differentiation with respect to $\tau$, and we introduced the notation
\begin{equation}
\alpha=\frac{k_0\tau_R^2}{m}.
\label{1vI(7.5)}
\end{equation}
Note that for the Maxwell model (see Section~\ref{1dsSection2}, Eq.~(\ref{1vI(2.2)})) we have
$k_0=k$, $\tau_R=b/k$, and $\alpha=1/(4\zeta^2)$.

Furthermore, according to Eqs.~(\ref{1vI(7.2)})), the variable impact velocity is
$$
\dot{x}(t)=v_0\xi^\prime(\tau).
$$

Let $\tau_c$ be the dimensionless duration of the impact process. Then, the coefficient of restitution can be found as
\begin{equation}
e_*=-\xi^\prime(\tau_c).
\label{1vI(7.6)}
\end{equation}

From Eqs.~(\ref{1vI(7.3)})) and (\ref{1vI(7.4)})), it is evident that $\tau_c$ is a function of $\alpha$ only and does not depend on $v_0$. Thus, in view of (\ref{1vI(7.6)})), we conclude that the coefficient of restitution $e_*$ is constant with respect to the initial impact velocity $v_0$. 

It can be shown that the same qualitative conclusions are drawn from the linear biphasic model \cite{MowKueiLaiArmstrong1980} for articular cartilage deformation. In this case, the parameter $\tau_R$, which enters Eqs.~(\ref{1vI(7.2)})), can be defined as a typical diffusion time 
$\tau_D=h^2/(\kappa H_A)$, where $h$ is the cartilage layer thickness, $\kappa$ is the cartilage permeability, and $H_A$ is the aggregate modulus.

{\rem{
Let us consider the question of applicability of the coefficient of restitution for diagnosis of the state of health of the tissue. In the framework of the asymptotic model (\ref{1vI(9.8)}), according to Eq.~(\ref{1vI(2.9)}), we will have
\begin{equation}
e_*=\exp\Bigl(-\frac{\pi\zeta}{\sqrt{1-\zeta^2}}\Bigr),
\label{1vI(e.1)}
\end{equation}
where (see the last formula (\ref{1vI(2.6)}))
\begin{equation}
\zeta=\frac{k}{2\omega_0 b}
\label{1vI(e.2)}
\end{equation}
with $k$ and $b$ defined as follows (see Eq.~(\ref{1vI(9.9)}) and (\ref{1vI(9.90)})):
\begin{equation}
k=\frac{3\mu_s a^4}{16h^3},\quad \frac{b}{k}=\tau_R=\frac{h^2}{3\mu_s \kappa}.
\label{1vI(e.3)}
\end{equation}
From (\ref{1vI(e.1)}) it is readily seen that the coefficient of restitution $e_*$ decreases with increasing loss factor $\zeta$. At the same time, in view of (\ref{1vI(e.2)}) and (\ref{1vI(e.3)}), we have
\begin{equation}
\zeta=\frac{\sqrt{3m\mu_s}\kappa}{2a^2\sqrt{h}}.
\label{1vI(e.4)}
\end{equation}
It is known \cite{Wu2000} that for articular cartilage in the early stages of osteoarthritis, the following degenerative changes are observed: increased permeability, $\kappa$, increased thickness of the cartilage layer, $h$, reduced shear modulus, $\mu_s$, and/or a combination of these effects.
Formula (\ref{1vI(e.4)}) implies that increasing the permeability of the cartilage results in a decrease of $e_*$, while increasing the cartilage thickness and decreasing the shear modulus (both create a softening effect) apparently results in an increase of $e_*$. Because it is known \cite{Knecht2006} that osteoarthritic cartilage may show a dramatic (more than 6-fold) increase in the hydraulic permeability $\kappa$, it can be expected that the overall change in the coefficient of restitution $e_*$ will be negative.

It should be also observed that formulas (\ref{1vI(e.4)}) and (\ref{1vI(2.9)}) imply an increase of the coefficient of restitution $e_*$ with increasing the cartilage thickness $h$, whereas the experimental data presented in Fig.~\ref{Varga2007.pdf} apparently show an inverse tendency.
}}

{\rem{
It is known \cite{VerteramoSeedhom2007,BurginAspden2008} that impact loading of articular cartilage at high impact stresses typically result in fissuring of the articular cartilage surface. At the same time, the formation of cracks allows to absorb greater amounts of energy as well as dramatically affect the deformation resistance of cartilage resulting in change of the parameters of the impact model. In other words, the mechanical properties of the tissue do not remain the same to the end of the impact process. 
}}

Observe that the biphasic theory incorporating Lam\'e parameters assumes that the material of solid phase is linearly elastic in order for these to have unique values. But if the material is viscoelastic these parameters are difficult to define and they become functions of deformation and/or time, if they are meaningful at all. Furthermore, there is an intrinsic circularity problem associated with using the aggregate modulus $H_A$, which is evaluated at equilibrium after the interstitial water is squeezed out, to define the mechanical properties that are then assumed to pertain during the impact deformation. 
Thus, the fact that the biphasic theory provides a good fit to measured curves in the creep and stress relaxation tests can be basically considered as a consequence of a curve-fitting procedure with a minimum of three free parameters rather than a derivation from first principles. In other words, it remains to be an open question on the efficiency of mixture models for articular cartilage at high strain rates. 

In the present study we addressed the question of whether the main features of articular impact observed in
\cite{Varga2007,Edelsten2010} could be qualitatively predicted using a linear viscoelastic theory or the linear biphasic theory.
It is to note that the deformations encountered in impact tests should be small enough for the linear theories to apply. With respect to engineering polymers note that the linear theory of viscoelasticity may hold reasonably well even up to some $5-10\%$ extension, in particular for certain rubbers 
\cite{Tschoegl1997}.

\section{Conclusions}
\label{1dsSectionC}

The results of this study based on the linear viscoelasticity imply the following properties of the linear impact models:

1. The coefficient of restitution $e_*$ is a function of the damping ratio $\zeta$ alone. This means that $e_*$ does not depend on the impact velocity $v_0$, but it depends on the impactor mass $m$ and the sample thickness (through the stiffness $k$).

2. The impact duration $t_c$ is inversely proportional to $\omega_0$, that is $t_c$ is proportional to $\sqrt{m}$, and depends on the damping ratio as well. The impact duration does not depend on the impact velocity $v_0$.

3. The maximum displacement, $x_m$, and the maximum contact force, $F_M$, are proportional to $\sqrt{m}$ and $v_0$.

4. The time variation of the incremental dynamic stiffness $dF/dx$ remains the same for different initial impact velocities. 

5. In the drop weight impact test, the gravitational effect increases the impact duration $t_c$ and decreases the coefficient of restitution $e_*$. At that, the coefficient of restitution increases with the impact velocity $v_0$.

\section*{Acknowledgment}

The financial support from the European Union Seventh Framework Programme under contract number PIIF-GA-2009-253055 is gratefully acknowledged. The author also would like to express his gratitude to the Referees for their helpful comments and discussions.


\end{document}